\documentclass[10pt]{article}
\usepackage{latexsym}
\usepackage{graphicx}
\usepackage{color}
\usepackage{amsmath}
\usepackage{amsfonts}
\usepackage{amssymb}
\usepackage{url}
\usepackage{wrapfig}

\definecolor{darkblue}{RGB}{0, 0, 102}
\definecolor{IKB}{RGB}{0, 47, 167}

\newtheorem{theorem}{Theorem}
\newtheorem{lemma}[theorem]{Lemma}

\newtheorem{cor}[theorem]{Corollary}
\newtheorem{example}[theorem]{Example}

\newtheorem{remark}[theorem]{Remark}

\newtheorem{Def}[theorem]{Definition}

\newcommand{\bmath}[1]{{\boldmath ${#1}$}}

\newcommand{\p}{\ensuremath{\mathcal{P}}}

\newcommand{\A}{\ensuremath{\mathcal{A}}}

\newcommand{\F}{\ensuremath{\mathcal{F}}}

\newcommand{\cL}{\ensuremath{\mathcal{L}}}

\newcommand{\V}{\ensuremath{\mathcal{V}}}
\newcommand{\I}{\ensuremath{\mathcal{I}}}

\newcommand{\cP}{\ensuremath{\mathcal{P}}}

\newcommand{\E}{\begin{equation}}
\newcommand{\EE}{\end{equation}}
\newcommand{\QED}{\ \rule{.1in}{.1in}}

\newcommand{\cZ}{\mathbb{Z}}
\newcommand{\R}{\mathbb{R}}

\newcommand{\conv}{\mathop{\rm conv}}

\advance \topmargin by -\headheight
\advance \topmargin by -\headsep
\begin{document}

\hspace{.1in}\\

\vskip 1cm

\begin{center}
\Large {\bf LP formulations for mixed-integer polynomial optimization problems}

\normalsize

Daniel Bienstock and Gonzalo Mu\~noz, Columbia University, December 2014
\end{center}

\begin{abstract}
We present a class of linear programming approximations for mixed-integer polynomial optimization problems that take advantage of structured sparsity of the 
constraint matrix. In particular, we show that if the intersection graph of the
constraints has tree-width bounded by a constant, then for any desired
tolerance there is a linear programming
formulation of polynomial size.  Via an additional reduction, 
 we obtain a polynomial-time 
approximation scheme for the ``AC-OPF'' problem on graphs
with bounded tree-width.  These constructions partly rely on a
general construction for pure binary optimization problems
where individual constraints are available through a membership oracle; if the intersection graph for the constraints has bounded tree-width our construction is of linear
size and exact.  This improves on a number of results in the literature, both from
the perspective of formulation size and generality.
\end{abstract}

\section{Introduction}
A fundamental paradigm in the solution of integer programming and combinatorial optimization 
problems is the use of extended, or lifted, formulations, which rely on the binary nature of
the variables and on the structure of the constraints to generate higher-dimensional convex
relaxations with provably strong attributes.  In this paper we consider mixed-integer polynomial optimization
problems.  We develop a reformulation operator
which relies on the combinatorial structure of the constraints to produce 
linear programming approximations which attain provable bounds.  A major
focus is on polynomial optimization problems over networks and our main
result in this context (Theorem
\ref{npotheorem_0} below) implies as a corollary that there exist polynomial-size linear
programs that approximate  the AC-OPF problem and the fixed-charge network flow
problem on bounded tree-width graphs. 

Our work relies on the concepts of \textit{intersection graph} and
\textit{tree-width}; as has been observed before (\cite{tw}, \cite{Laurentima}, \cite{lasserre2}, \cite{Waki06sumsof}, \cite{wainjor}), the combination of these two concepts makes it possible to define a notion of structured sparsity in an 
optimization context that we will exploit here (see below for more references). The intersection graph of a system of constraints
is a central concept originally introduced in \cite{fulkersongross} and which has been
used by many authors, sometimes using different terminology.

\begin{Def} \label{intersectgraph} The {\bf intersection graph} of a system of constraints is the undirected
graph which has a vertex for each variable and an edge for each pair of variables 
that appear in any common constraint.  
\end{Def}
\begin{example}\label{intersectex} Consider the system of constraints on
variables $x_1, x_2, x_3$ and $x_4$.
\begin{equation} \label{iex}
\begin{array}{lll}
 (1) \ 3 x_1^2  - \log_2(1 + x_2)  \ge 0, \quad & (2) \ -2 x_2^2 \ + \ (1+ x_3)^3 \ \ge \ 0, \quad & (3)\ x_3x_5 = 0 \\
 (5) \ x_3^3 \ - \ 2 - x_4^2 \ < \ 0, \quad & (5) \ (x_1, x_4) \in A &
\end{array}
\end{equation}
Where $A$ is some arbitrary set. Then the intersection graph has vertices $\{1,2,3,4,5\}$ and edges $\{1,2\}$, 
$\{2,3\}$, $\{3,4\}$, $\{3,5\}$ and $\{4,1\}$. 
\end{example}

\begin{Def} \label{tree-widthdef} An undirected graph has {\bf tree-width} $\le k$ if it is contained in a chordal graph with clique number $\le k+1$.
\end{Def}
The tree-width concept was explicitly defined in \cite{Robertson198449} (also see \cite{rs86}), but there are many equivalent definitions.  An earlier discussion is found in \cite{halin} and closely related concepts have been used by many authors under other names, e.g. the ``running intersection'' property, and the notion of ``partial k-trees".  An important known fact is that an $n$-vertex graph with tree-width $\leq k$ has $O(k^2 n)$ edges, and thus low tree-width graphs are \emph{sparse}, although the converse is not true. In the context of this paper, we can exploit structural sparsity in 
an optimization problem when the tree-width of the intersection graph is small. 
We also note that bounded tree-width can 
be recognized in linear time \cite{bodlaender2}.  

Our focus, throughout, is on obtaining polynomial-size LP formulations.  The
construction of such ``compact'' formulations is a goal of fundamental 
theoretical importance and quite separate from the development of polynomial-time algorithms.  From
a numerical perspective, additionally, the representation of a problem as an LP
permits the use of practical bounding techniques such as cutting-plane 
and column-generation algorithms. We first prove:
\begin{theorem} \label{potheorem_0}  Consider a mixed-integer, linear objective, polynomially constrained problem
\begin{subequations} \label{PO}
\begin{eqnarray}
\mbox{(PO):} && \min ~ c^T x \\
\mbox{subject to:} && \quad \quad f_i(x) \ \ge \ 0 \quad 1 \le i \le m\\
&& \quad \quad x_j \in \{0,1\} \quad 1 \le j \le p, \quad x_j \in [0,1] \quad p+1 \le j \le n. 
\end{eqnarray}
\end{subequations}
Let $\pi$ and $F$ be such that for $0 \le i \le m$,  $f_i(x)$ has maximum degree at most $\pi$ and 
$L_1$-norm of coefficients at most $F$. If
the intersection graph of the constraints has tree-width $\le \omega$ then for
any  $0 < \epsilon < 1$, there is a linear programming formulation with $O\left((2 \pi/\epsilon)^{\omega + 1} \, n \log (\pi/\epsilon)\right)$ variables and constraints that solves $\cP$ within feasibility tolerance $F \epsilon$ and optimality tolerance $\|c\|_1 \epsilon$.
\end{theorem}

Below we will provide an extended statement for this result, as well as a precise
definition of `tolerance'. However, the statement in Theorem \ref{potheorem_0} is
indicative of the fact that as $\epsilon \rightarrow 0$ we converge to an optimal solution, and the computational workload
grows proportional to $O(\epsilon^{-\omega - 1} \log \epsilon^{-1})$.  Moreover,
as we will argue, it is straightforward to prove that unless $P=NP$, no polynomial time algorithm for mixed-integer polynomial optimization exists that improves on the dependence on $\epsilon$ given by Theorem \ref{potheorem_0}.
As far as we know
this theorem is the first to provide a polynomial-size formulation 
for polynomial optimization problems with guaranteed bounds.  

Our next result is motivated by recent work on the AC-OPF (Optimal Power Flow) problem in electrical transmission
\cite{lavaeilow2012}, \cite{boselowqcqp12}, \cite{javadsojoudi14}.  A generic version of this problem
can be succinctly described as follows.  We are given an undirected graph
$G$ where for each vertex $u \in V(G)$ we have two
variables, $e_u$ and $f_u$.  Further, for each edge $\{u,v\}$ we have four
$2 \times 2$ matrices $M_{uv}, M_{vu}, N_{uv}$ and $N_{vu}$.  For 
each edge $h = \{u,v\}$, write $w_h = (e_u, f_u, e_v, f_v)^T$. Then we have
\begin{subequations} \label{acopf}
\begin{eqnarray}
\mbox{(AC-OPF):} && \min \quad \sum_{u \in V} c_u g^2_u \\
\mbox{subject to:} && L_u \ \le \ \sum_{h = \{u,v\} \in \delta(u)} w_h^T M_{uv} w_h \ \le \ U_u \quad \forall u \in V(G) \label{active}\\
&& L'_u \ \le \ \sum_{h = \{u,v\} \in \delta(u)} w_h^T N_{uv} w_h \ \le \ U'_u \quad \forall u \in V(G) \label{reactive}\\
&& B_u \ \le \ ||(e_u, f_u)^T|| \ \le \ A_u \quad \forall u \in V(G)\\
&& g_u \ = \ \sum_{h = \{u,v\} \in \delta(u)} w_h^T M_{uv} w_h \quad \forall u \in V(G). \label{acopflast}
\end{eqnarray}
\end{subequations}
In this formulation $L_u, L'_u, U_u, U'_u, B_u, A_u$ and $c_u$ are given values,
and the $g_u$ are auxiliary variables defined as per (\ref{acopflast}). Here and
below, given a graph $H$ its vertex set is $V(H)$ and its edge set is $E(H)$,
and for 
$u \in V(H)$ we use $\delta(u) = \delta_H(u) $ to denote the set of edges incident with $u$ and write $\deg(u) = \deg_H(u) = | \delta(u) |$.

As a generalization of AC-OPF, we consider \textit{network mixed-integer polynomial
  optimization problems} ({\bf NPO}s, for short).  These are
PO problems with an underlying network structure
specified by a graph $G$.  Specifically, we assume for each
vertex $v$ in $G$ there is a set $X_v$ of variables \textit{associated}
with $v$. Moreover,
each constraint is associated with one vertex of 
$G$;  a constraint associated with vertex $u$ takes the form 
\begin{eqnarray}  \sum_{\{u,v\} \in \delta(u)} p_{u,v}(X_{u} \cup X_v) & \geq & 0, \label{networkconstraint}
\end{eqnarray}
where $\delta(u)$ is as defined above and each $p_{u,v}$ is 
a polynomial. Note that this definition allows a vertex $v$ to have many constraints of the type \eqref{networkconstraint} associated with it. The sets
$X_v$ are not assumed to be pairwise disjoint, and thus a given variable may appear in several
such sets; however, for technical reasons, we assume that for any variable $x_j$ the set $\{ v \in V(G) \, : \, x_j \in X_v\}$ induces a connected subgraph of $G$. Clearly AC-OPF is an NPO (with the $X_v$ pairwise disjoint),  and it can also be shown that
optimization problems on gas networks \cite{lars} are NPOs, as well.

Yet 
another example is provided by the classical
\textit{capacitated fixed-charge network flow} problem (see \cite{fcnf}) which
has received wide attention in the mixed-integer programming literature.  In
the simplest case we have a directed graph $G$; for each vertex $u$ we have a
value $b_u$ and for each arc $(u,v)$ we have values $f_{uv}$, $c_{uv}$ and $w_{uv}$. The problem is
\begin{subequations} \label{fcnf}
\begin{eqnarray}
\mbox{(FCNF):} && \min \quad \sum_{(u,v)} f_{uv} y_{uv} \, + \, c_{uv} x_{uv}  \\
\mbox{subject to:} && N x \ = \ b, \\
&& 0 \, \le \, x_{uv} \, \le \, w_{uv} y_{uv}, \quad y_{u,v} \in \{0,1\} \quad \forall (u,v).
\end{eqnarray}
\end{subequations}
where $N$ is the node-arc incidence matrix of $G$. This is an NPO (e.g. associate each $x_{u,v}$ with either $u$ or $v$). When $G$ is a \textit{caterpillar} (a  path with pendant edges) it includes the knapsack problem as a special case.  FCNF can arise
in supply-chain applications, where $G$ will be quite sparse and often tree-like.

\begin{figure}[h] 
\centering
\includegraphics[height=1in]{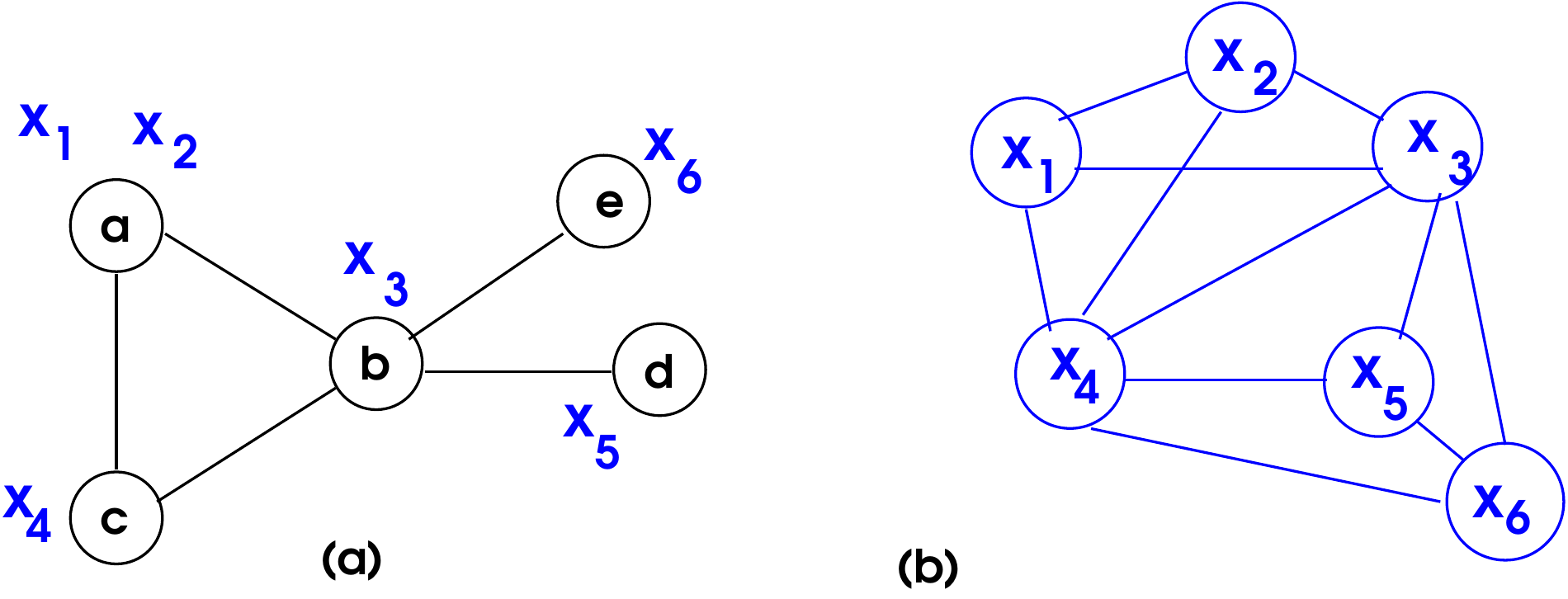}
\caption{(Example \ref{networkex}) (a) Network polynomial optimization problem. (b) The intersection graph.\label{fig:network}}
\end{figure}

Above (Theorem \ref{potheorem_0}) we have focused on exploiting the 
structure of the intersection graph for a problem; as we discussed this graph is obtained \textit{from a formulation} for the problem. However, in NPOs there is 
already a graph, which in the above examples frequently has
moderate tree-width, and it is this condition that we would like to exploit.
In fact, recent work (\cite{molzahnsdpimplementation}, \cite{javadimplementation}) develops faster solutions to SDP \textit{relaxations} of AC-OPF problems by leveraging small tree-width of the underlying graph. To highlight the difference between the two
graphs, consider the following examples:
\begin{example}\label{networkex}  Consider the NPO with constraints
\begin{eqnarray}
&& x_1^2 + x_2^2 + 2 x_3^2 \ \le \ 1, \quad x_1^2 - x_3^2 + x_4\ \ge \ 0, \nonumber \quad x_3 x_4 + x_5^3 - x_6\ \ge \ 1/2 \nonumber \\
&& 0 \le x_1 \le 1,  \quad 0 \le x_2 \le 1, \quad 0 \le x_3 \le 1, \quad 0 \le x_5 \le 1, \quad 0 \le x_6 \le 1, \quad x_4 \in \{0, 1\}. \nonumber
\end{eqnarray}
Here, variables $x_1$ and $x_2$ are associated with vertex $a$ of the graph
in Figure \ref{fig:network}, and variables $x_3$, $x_4$, $x_5$ and $x_6$ are associated with vertices $b, c, d$ and $e$, respectively.  The first constraint
involves variables associated with the endpoints of edge $\{a,b\}$, 
the second concerns edges $\{a,b\}$ and $\{a, c\}$, and the third concerns edges
$\{b,c\}, \{b,e\}$ and $\{b, d\}$.
\end{example}
\begin{example}\label{minknapsackex} Consider a knapsack problem 
$\min \{ \, c^Tx \ : \ a^Tx \ \ge \ b, \ x \in \{0,1\}^n \, \}$.  This is a NPO, using a star network on
$n+1$ vertices, which has tree-width $1$.  Yet, 
  if $a_j \neq 0$ for all $j$, the intersection graph is a clique of size $n$.  Note that we can restate $x \in \{0,1\}^n$ as $x_j(1-x_j) = 0, \, \forall j$.
\end{example}
In fact, even an AC-OPF instance, on a tree, can give rise to an intersection graph
with high tree-width, because constraint (\ref{active}) or (\ref{reactive}) mirrors the intersection graph
behavior in Example \ref{minknapsackex}.\\

In the general case, in order to avoid the possible increase in tree-width going from the underlying graph to the intersection graph, we show that given a NPO problem on a graph of small tree-width, there is an equivalent 
NPO problem whose intersection graph also has small tree-width.  As a result of
this elaboration we obtain:

\begin{theorem} \label{npotheorem_0} Consider a network mixed-integer polynomial
optimization problem over a graph $G$ of tree-width $\le \omega$ and maximum-degree $D$, over $n$ variables, and where every polynomial $p_{u,v}$ has maximum degree $\le \pi$.  Suppose that the number of variables plus the number of constraints associated with any vertex of $G$ is at most $\Delta$.  Given $0 < \epsilon < 1$, there is a linear programming formulation of size $O((D \pi/\epsilon)^{O(\Delta \omega)}\, n \, \log (\pi/\epsilon))$ that solves the problem within scaled tolerance $\epsilon$.
\end{theorem}

\noindent Here, ``scaled tolerance'' embodies the same notion of optimality and feasibility approximation as in 
Theorem \ref{potheorem_0}.  Later we will discuss why this approximation
feature is needed. Further, we note that in the case of the AC-OPF problem we
have $\Delta = 5$ and $\pi = 2$.

\begin{cor} There exist polynomial-time
approximation schemes for the AC-OPF problem on graphs with bounded tree-width,
and for the capacitated fixed-charge network flow problem on
graphs with bounded tree-width.
\end{cor}

Our third result is important toward the proof of Theorem \ref{potheorem_0}
but is of independent interest. As we will see, the construction in 
Theorem \ref{potheorem_0}
approximates a mixed-integer polynomially constrained problem  with a polynomially constrained pure binary polynomially constrained problem.  As a generalization, 
we study ``general'' binary
problems, or {\bf GB} for short, defined as follows.
\begin{itemize}
\item [(i)] There are $n$ variables and $m$ constraints.  For $1 \le i \le m$, constraint $i$ is 
characterized by a subset $K[i] \subseteq \{1,\ldots, n\}$ and a set $S^i \subseteq \{0,1\}^{K[i]}$.
Set $S^i$ is implicitly given by a \textit{membership oracle}, that 
is to say a mechanism that upon input $y \in \{0,1\}^{K[i]}$, truthfully reports
whether $y \in S^i$.
\item [(ii)] The problem is to minimize a linear function $c^Tx$, over $x \in \{0,1\}^{n}$, and subject to the
constraint that for $1 \le i \le m$ the sub-vector $x_{K[i]}$ is contained in $S^i$.
\end{itemize}
Any linear-objective, binary optimization problem  whose
constraints are explicitly 
stated can be recast in the form GB; e.g., each set $S^i$ could be
described by a system of algebraic equations in the variables $x_j$ for $j \in K[i]$.  GB problems are related to classical constraint satisfaction problems, however the
terminology above will prove useful later.  A proof of part (a) in the following
result can be obtained using techniques in \cite{Laurentima} (Section 8)
although not explicitly stated there.  We will outline this proof,
which relies on the ``cone of set-functions'' approach of \cite{Lovasz91conesof}
and
also present a new proof.

\begin{theorem} \label{genbtheorem_0} Consider a
GB problem whose intersection graph has tree-width $\le \omega$. 
\begin{itemize}
\item[(a)] There is an exact linear programming formulation with $O(2^\omega n)$ variables and constraints, with $\{0,1,-1\}$-valued constraint coefficients. 
\item[(b)] The formulation 
can be constructed by performing $2^{\omega}m$ oracle queries and with 
additional workload $O^*(\omega 2^{\omega} \sum_{t \in V(T)} |\left\{i \, : \, K[i] \subseteq Q_t \right\}| + \omega mn) = O^*(\omega 2^{\omega}mn + \omega^2 2^{\omega}n)$, where the ``*'' notation indicates logarithmic factors in $m$ or $n$.
\end{itemize}
\end{theorem}

 Note that the size of the formulation is independent of the number constraints in the given instance of GB.   And even though we use the general 
setting of membership oracles, this theorem gives an exact reformulation, as opposed to Theorems \ref{potheorem_0} and \ref{npotheorem_0}, where an approximation is required
unless P=NP. Theorem \ref{genbtheorem_0} has additional implications toward linear and polynomial binary optimization problems. We will examine these issues in Section \ref{purebinary}. Regarding part (b) of the theorem, it can be
shown that $2^{\omega}m$ is a \textit{lower bound} on the number of oracle
queries that any algorithm for solving GB must perform.\\

Theorem \ref{genbtheorem_0} describes an \textit{LP formulation} for GB; the
relevance of this focus was discussed above.
Together with the reductions used to obtain Theorems \ref{potheorem_0} and 
\ref{npotheorem_0} we obtain approximate LP formulations for polynomial 
(resp., network polynomial) mixed-integer problems.  
Of course, Theorem \ref{genbtheorem_0} also
implies the existence of an \textit{algorithm} for solving GB in time polynomial in
$O(2^\omega (n + m))$.  However one
can also derive a direct algorithm of similar complexity using well-known,
prior ideas on polynomial-time methods for combinatorial 
problems on graphs of bounded tree-width.

\subsubsection{Prior work}\label{prior}
There is a broad literature dating from the 1980s on polynomial-time
algorithms for combinatorial
problems on graphs with bounded tree-width.  
An early reference is \cite{arnborg0}. Also see \cite{arnborg}, \cite{acp}, \cite{BrownFellowsLangston1989}, 
\cite{Bern1987216}, \cite{bodlaender}, \cite{bienstocklangston} and
from a very general perspective, \cite{tovey}.
These algorithms rely on ``nonserial dynamic programming'', i.e., dynamic-programming on trees. See \cite{george2}, \cite{george}, \cite{nonserial}.

A parallel research stream concerns ``constraint satisfaction problems'', or CSPs. Effectively, the feasibility version of problem GB is a CSP.  One
can also obtain an algorithm for problem GB, with similar complexity,  and
relying on similar dynamic programming ideas as the algorithms above, 
from the perspective of \textit{belief propagation} on an appropriately defined \textit{graphical model}. Another central technique
is the tree-junction theorem of \cite{graphical}, which shows how a 
a set of marginal probability distributions on the edges of a hypertree 
can be extended to a joint distribution over the entire vertex set. Early references are \cite{pearl82}, \cite{freuder} and \cite{dechterpearl}. Also see \cite{wainjor}, \cite{chandra} (and references therein), and \cite{wainwrightjordanlong}.

Turning to the integer programming context, \cite{bz} (also see the PhD thesis \cite{zuckerberg}) develop extended formulations for binary linear programs by considering the subset algebra of 
feasible solutions for individual constraints or small groups of constraints;
this entails a refinement of the 
cone of set-functions approach of \cite{Lovasz91conesof}. The method in \cite{bz} is similar to the one in this paper, in that here
we rely on a similar algebra and on extended, or ``lifted'' reformulations for $0/1$ integer programs.  The classical examples in this vein 
are the reformulation-linearization technique of \cite{sheraliadams}, the cones of matrices method \cite{Lovasz91conesof}, the lift-and-project method of  \cite{balasceriacornuejols}, and the moment relaxation methodology of \cite{lasserre}.  See \cite{Laurent01acomparison} for a unifying analysis; another comparison is
provided in \cite{AuTuncel13}.

The work in \cite{tw} considers \textit{packing binary integer programs} are considered, i.e.
problems of the form 
\begin{eqnarray}
&& \max\{c^T x \, : \, Ax \le b, \ x \in \{0,1\}^n\} \label{packing}
\end{eqnarray}
where $A \ge 0$ and integral and $b$ is integral. Given a valid inequality $\alpha x \ge \beta$, its
\textit{associated graph} is the subgraph of the intersection graph induced
by $S = \{ 1 \le j \le n: \, \alpha_j \neq 0\}$;  i.e. it has vertex-set
$S$ and there an
edge $\{j, k\} \in S \times S$ whenever $a_{ij} \neq 0$ and $a_{ik} \neq 0$ for some row $i$.

In \cite{tw} it is shown that given and $\omega \ge 1$, the level-$\omega$ Sherali-Adams reformulation of (\ref{packing})  implies every valid inequality whose associated graph has tree-width $\le \omega - 1$.  Further, if $A$ is $0/1$-valued,
the same property holds when the associated graph has tree-width $\le \omega$.
As a corollary, given a graph $G$ with tree-width $\le \omega$, the Sherali-Adams reformulation of the vertex 
packing linear program $\{ x \in [0,1]^{V(G)} \, : \, x_u + x_v \le 1 \ \forall \, \{u,v\} \in E(G)\}$, which has $O(\omega n^{\omega +2})$ variables and constraints, is exact.  As far as we know this is the first result linking 
tree-width and reformulations for integer programs.   A different
result which nevertheless appears related is obtained in \cite{cunngeel}.  

In \cite{wainjor}, binary polynomial optimization problems are considered, i.e problems as
$ \min \{ \, c^T x \, : \, x \in \{0,1\}^n, \ g_i(x) \ge 0, \ 1 \le i \le M\}$ 
where each $g_i(x)$ is a polynomial. They show that if the tree-width of the 
intersection graph of the constraints is $\le \omega$, then the  level-$\omega$ Sherali-Adams or Lasserre reformulation of the problem is exact.
Hence there is an LP formulation with $O(n^{\omega + 2})$ variables and $O(n^{\omega + 2} M)$ constraints.

A comprehensive survey of results on polynomial optimization and related
topics is provided in \cite{Laurentima}.  Section 8 of \cite{Laurentima} builds on the work in \cite{Laurent01acomparison}, which provides a common framework
for the Sherali-Adams,  Lov\'{a}sz-Schrijver and Lasserre reformulation operators. In addition to the aforementioned results related to problem GB (to which
we will return later), \cite{Laurentima} explicitly shows that the 
special case of the vertex-packing
problem on a graph with $n$ vertices and  tree-width $\le \omega$ has a formulation of size $O(2^{\omega} n)$; this is stronger than the implication from \cite{tw} discussed above.  Similarly, it is shown in \cite{Laurentima} that the
max-cut problem on a graph with $n$ vertices and  tree-width $\le \omega$ has a formulation of size $O(2^{\omega} n)$.

In the continuous variable polynomial optimization setting, \cite{kojima2}, and \cite{Waki06sumsof}
present methods for exploiting low tree-width of the intersection matrix
e.g. to speed-up the sum-of-squares or moment relaxations of a problem. Also see \cite{Grimm} and Section 8 of \cite{Laurentima}.
\cite{lasserre2} considers polynomial optimization problems as well. 
In abbreviated form, \cite{lasserre2} shows that where $p$
is the tree-width of the intersection graph, there is a hierarchy of semidefinite relaxations
 where the $r^{th}$ relaxation ($r = 1, 2, \ldots ...$) has $O( n p^{2r} )$ variables and $O(n + m)$ LMI constraints; further, as $r \rightarrow +\infty$ the value of the
relaxation converges to the optimum. Also see \cite{Nie2014} and \cite{lasserre2015}.

Finally, there are a number of results on using lifted formulations for polynomial
optimization problems, along the lines of the RLT methodology of \cite{sheraliadams}.
See \cite{sheralitunc}, \cite{sheralidalk} and references therein.

\subsubsection{Organization of the paper}
This paper is organized as follows.  Mixed-integer polynomial optimization problems and a proof of Theorem \ref{potheorem_0} 
are covered in Section \ref{minpoly}.  Network mixed-integer polynomial
optimization problems and Theorem \ref{npotheorem_0} are addressed in 
Section \ref{graphminpoly}. And finally,  in Section \ref{purebinary}, we will present a detailed analysis of the pure
binary problems addressed by Theorem \ref{genbtheorem_0} and a proof of this result.

\section{Mixed-integer polynomial optimization problems}\label{minpoly}

In this section we consider mixed-integer polynomial optimization problems (PO) and 
prove a result (Theorem \ref{potheorem}, below) that directly implies Theorem \ref{potheorem_0} given in the introduction. This proof will make use of the result
in Theorem \ref{genbtheorem_0}, to be proven in Section \ref{purebinary}.

In what follows we will rely on a  definition of tree-width which is equivalent to Definition \ref{tree-widthdef}. This definition makes use of the concept of \emph{tree-decomposition}:

\begin{Def} Let $G$ be an undirected graph.  A {\bf tree-decomposition} \cite{Robertson198449}, \cite{rs86} of $G$ is a pair $(T, Q)$ 
where $T$ is a tree and $Q = \{ Q_t \, : \, t \in V(T) \}$ is a family of
subsets of $V(G)$ such that
\begin{itemize}
\item[(i)] For all $v \in V(G)$, the set $\{ t \in V(T) \, : \, v \in Q_t \}$
forms a subtree $T_v$ of $T$, and
\item[(ii)] For each $\{u, v \} \in E(G)$ there is a $t \in V(T)$ such
that $\{u, v \} \subseteq Q_t$, i.e. $t \in T_u \cap T_v$.
\end{itemize}
The {\em width} of the decomposition is $\max \left\{ | Q_t | \, : \, t \in V(T) \right\} \, - \, 1$.  The tree-width of $G$ is the minimum width of a 
tree-decomposition of $G$.  See Example \ref{tree-widthex}.
\end{Def}
 
\begin{figure}[h] 
\centering
\includegraphics[width=5in]{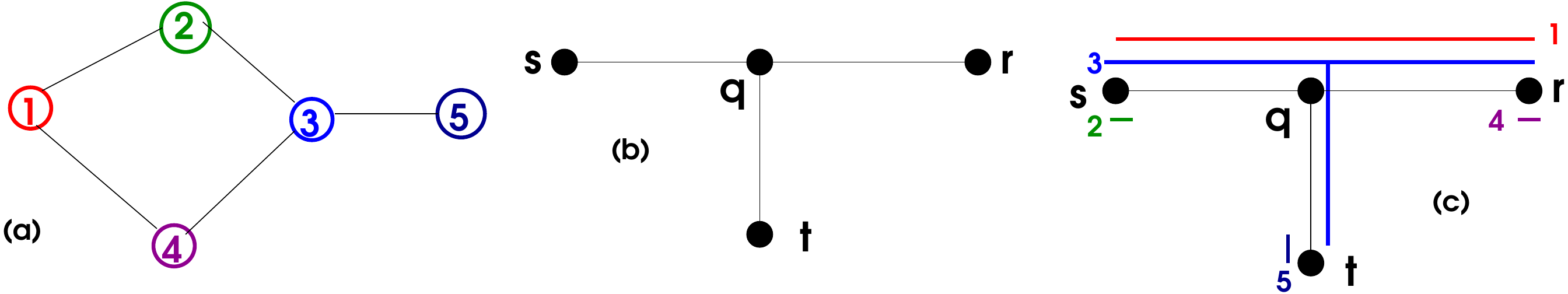}
\caption{A tree-decomposition. (a) Graph $G$. (b) Tree $T$.  (c) Tree $T$ with subtrees $T_v$.\label{fig:treedecomp}}
\end{figure}

\begin{example}\label{tree-widthex} {\bf (Tree-decomposition)}
Consider the intersection graph $G$ arising in Example \ref{intersectex}. See Figure \ref{fig:treedecomp}(a).  A tree-decomposition with tree
$T$ is shown in Figure \ref{fig:treedecomp}(b)-(c).
\end{example}

Since this definition relates a specific decomposition of the graph with its tree-width, many of the arguments we provide will rely on modifying or creating valid tree-decompositions that attain the desired widths.

We make some remarks pertaining to to the PO problem.  Throughout we will use the definition of
used in the introduction (formulation (\ref{PO})). The $i^{th}$
constraint, for $1 \le i \le m$, is given by $f_i(x) \geq 0$, where $f_i(x)$ is a has the form
\begin{eqnarray}
  f_i(x) & = & \sum_{\alpha \in I(i)} f_{i,\alpha} x^\alpha. \label{genpoly2}
\end{eqnarray}
Here $I(i)\subseteq \{1, \ldots, n\}$ is a finite set, $f_{i,\alpha}$ is rational and $x^\alpha$ is a monomial in $x$:
$$ x^\alpha \ = \ \prod_{j=1}^n x_j^{\alpha_j}.$$
Finally, for each $1\leq i \leq m$ we will denote as $\|f_i\|_1$ the 1-norm of the coefficients of polynomial $f_i(x)$, i.e
$$\|f_i\|_1  = \sum_{\alpha \in I(i)} |f_{i,\alpha}|.$$

Any linear-objective mixed-integer polynomial optimization problem where the feasible region is compact can be reduced to the
form (\ref{PO}) by appropriately translating and scaling variables.  

\begin{remark}\label{linearpo2} A polynomial-optimization problem with nonlinear objective can trivially be made into the form (\ref{PO}), by (for example) 
using a new variable and two constraints to represent each monomial in the objective.  Of course, such a modification
may increase the tree-width of the intersection graph.
\end{remark}

Now we precisely define what the intersection graph would be in this context.

\begin{Def} Given an instance of problem PO, let its intersection graph $\Gamma$ be the undirected
graph with $n$ vertices and where for $1 \le i \le m$ the set
$$\{ j\ : \ \alpha_j \neq 0,\ f_{i,\alpha} \neq 0,\ \alpha \in I(i) \}$$ 
induces a clique.  
\end{Def}
\begin{Def} \label{podefstuff}Consider an instance of problem PO. 
\begin{itemize}
\item[(a)] Given $\epsilon > 0$, we say a vector $x \in \{0,1\}^{p} \times [0,1]^{n-p}$ is {\bf scaled-\bmath{\epsilon} feasible} if
\begin{eqnarray}
f_i(x) & \ge & - \epsilon \|f_i\|_1, \quad 1 \le i \le m. \label{scaledfeasdef}
\end{eqnarray}
\item[(b)] We set \bmath{\pi \doteq \max_{1\leq i \leq m} \max_{ \alpha \in I(i)} \sum_{j=p+1}^n \alpha_j}. 
\end{itemize}
\end{Def}
\noindent We will prove the following result:
\begin{theorem} \label{potheorem} Given an instance of PO, let $\omega$ be the width of a tree-decomposition of the intersection graph.  For every $0 < \epsilon < 1$ there is a linear program 
$$LP \, : \, \min\{\hat c^T y \, : \, \hat A y \ge \hat b\}$$ 
with the following properties:
\begin{itemize}
\item [(a)] The number of variables and constraints is $O\left((2 \pi/\epsilon)^{\omega + 1} n \log (\pi/\epsilon)\right)$, and all coefficients are of polynomial size.
\item [(b)] Given any feasible solution $x$ to PO, there is a feasible solution
$y$ to LP with 
$$\hat c^T y \le c^T x + \epsilon \|c\|_1.$$
\item [(c)] Given an optimal solution $y^*$ to LP, we can construct
$x^* \in \{0,1\}^p \times [0,1]^{n-p}$ such that:
\begin{eqnarray}
\mbox{1.} & & \mbox{$x^*$ is scaled-$\epsilon$ feasible for PO, and} \nonumber\\
\mbox{2.} && c^T x^* \ = \ c^T y^*.\nonumber
\end{eqnarray}
\end{itemize}
\end{theorem}
\begin{remark} Assume PO is feasible. Then (b) shows LP is feasible and
furthermore by (c) solving LP yields a near-feasible solution for PO which 
may be superoptimal, but is not highly suboptimal.  Condition (a) states that the formulation $LP$ is of pseudo-polynomial size.
\end{remark}

To prove Theorem \ref{potheorem} we will rely on a technique used in \cite{glover75}; also see \cite{histo} and \cite{dash2007mir},  \cite{deybilinear} and citations therein.  
Suppose that $0 \le r \le 1$.  Then we can approximate $r$ as a sum of inverse
powers of $2$. Let $0 < \gamma < 1$ and
$$L \ = \ L(\gamma) \ \doteq \ \lceil \log_2 \gamma^{-1}\rceil.$$
\noindent Then
there exist $0/1$-values $z_h$, $1 \le h \le L$, with
\begin{eqnarray}\label{binapprox}
&& \sum_{h = 1}^L 2^{-h}z_h \ \le \quad r \quad \le \ \sum_{h = 1}^L 2^{-h}z_h \, + \, 2^{-L} \ \le \sum_{h = 1}^L 2^{-h}z_h \, + \, \gamma \, \le \, 1.  \label{sumof2}
\end{eqnarray}
Next we approximate problem PO with a problem of type GB.   For each $1 \le i \le m$ and  $\alpha \in I(i)$ we write
$$  Z(i,\alpha) = \{ j \ : \alpha_j \neq 0, \ 1\leq j \leq p \}.$$
In other words, set $Z(i,\alpha)$ is the set given by the indices of the \emph{binary} variables for PO that appear explicitly in monomial $x^\alpha$; thus for $j\in Z(i,\alpha)$ we have $x_j^{\alpha_j} = x_j$. Let 
\begin{eqnarray}
&& \mbox{ \bmath{\delta \ = \ \delta(\gamma, \pi) \ \doteq \ 1 - (1 - \gamma)^{\pi}}}. \label{deltadef}
\end{eqnarray}
The theorem will be obtained by using the following formulation, for appropriate
$\gamma = \gamma(\epsilon):$
\begin{subequations} \label{POgamma}
\begin{eqnarray}
\hspace{-.1in} \mbox{(GB$(\gamma)$)}: && \min ~ \sum_{j =1}^p c_j x_j \ + \ \sum_{j =p+1}^n c_j \left( \sum_{h = 1}^L 2^{-h} z_{j,h}\right)  \nonumber \\
\hspace{-.1in} \mbox{s.t.} && \sum_{\alpha \in I(i)} f_{i,\alpha} \left[ \, \prod_{j \in Z(i,\alpha)} x_j \, \prod_{j = p+1}^n  \left( \sum_{h = 1}^L 2^{-h} z_{j,h} \right)^{\alpha_j} \, \right] \ \ge \ - \delta \, \|f_i\|_1, \quad 1 \le i \le m \nonumber \\
\hspace{-.1in} && x_j \, \in \, \{0,1\}, \ 1 \le j \le p, \quad z_{j,h} \, \in \, \{0,1\}, \ \forall j \in \{p+1,\ldots, n\} , \ 1 \le h \le L.\nonumber
\end{eqnarray}
\end{subequations}
{\bf Remark.} This formulation replaces, in PO, each continuous variable $x_j$ with a sum of powers of $1/2$, using the binary variables $z_{j,h}$ in order to
effect the approximation (\ref{sumof2}).    \\

To prove the desired result we first need a technical property.
\begin{lemma}\label{tech} Suppose that for $1 \le k \le r$ we have values 
$u_k\geq 0, \ v_k \geq 0, \ q_k \in \cZ_+$ with $u_k + v_k \le 1$.  Then
$$ \prod_{k = 1}^r (u_k + v_k) ^{q_k} \, - \, \prod_{k = 1}^r u_k^{q_k} \ \le \ 1 \, - \,\prod_{k = 1}^r (1 - v_k)^{q_k}.$$
\end{lemma}
\noindent {\em Proof.} Take any fixed index $1 \le i \le r$. The expression
$$ \prod_{k = 1}^r (u_k + v_k) ^{q_k} - \prod_{k = 1}^r (u_k)^{q_k}$$
is a nondecreasing function of $u_i$ when all $u_k$ and $v_k$ are nonnegative, and so in the range $0 \le u_i \le 1 - v_i$ it is maximized when $u_i = 1 - v_i$. \QED\\

Using this fact, we can now show:
\begin{lemma}\label{tobinary} \hspace{.1in}
\begin{itemize}
\item[(a)] Suppose $\tilde x$ is a feasible for PO.  Then there is 
feasible solution for GB$(\gamma)$ with objective value at most $c^T \tilde x \ + \ \
\delta \|c\|_1$.
\item[(b)] Conversely, suppose $(\hat x, \hat z)$ is feasible for GB$(\gamma)$. Writing, for each $p+1 \le j \le n$,
$\hat x_j = \sum_{h = 1}^L 2^{-h} \hat z_{j,h}$, we have 
that $\hat x$ is scaled-$\delta$-feasible for PO and 
$c^T \hat x \, = \, \sum_{j =1}^n c_j \left( \sum_{h = 1}^L 2^{-h} \hat z_{j,h}\right)$.
\end{itemize}
\end{lemma}
\noindent {\em Proof.} (a) For each $j$ choose binary values $\tilde z_{j,h}$ so as to attain the approximation in (\ref{sumof2}).  Then for each $1 \le i \le m$ and $\alpha \in I(i)$ we have
\begin{eqnarray}
\hspace{-25pt} && \prod_{j = p+1}^n  \left( \sum_{h = 1}^L 2^{-h} \tilde z_{j,h} \right)^{\alpha_j} \ \le \ 
\prod_{j = p+1}^n  \tilde x_j^{\alpha_j} \ \le \ \prod_{j =p+1}^n  \left( \sum_{h = 1}^L 2^{-h} \tilde z_{j,h} \right)^{\alpha_j} \ + \ \delta. \nonumber
\end{eqnarray}
Here the left-hand inequality is clear, and the right-hand inequality follows
from Lemma \ref{tech} and the definition (\ref{deltadef}) of $\delta$.
Thus $\tilde z$ is feasible for $GB(\gamma)$ and the second assertion is
similarly proved.\\
\noindent (b) Follows by construction. \QED

We can now complete the proof of Theorem \ref{potheorem}.  Given an instance of 
problem PO together with a tree-decomposition of its intersection graph, of width $\omega$, we consider formulation GB$(\gamma)$ for $\gamma = \epsilon \pi^{-1}$.  
As an instance of GB, the formulation has at most $n L(\gamma)$ variables and its intersection graph has width at most $(\omega + 1)L(\gamma) - 1$.  To see this point, consider a
tree-decomposition $(T, Q)$ of the intersection graph for PO.  Then we obtain
a tree-decomposition $(T, \breve Q)$ for GB$(\gamma)$ by setting, for each $t \in V(T)$,
$$ \breve Q_t \ = \ \{ z_{j,h} \, : \, 1 \le h \le L(\gamma), \, j \in Q_t \}.$$

We then apply, to this instance of GB, Theorem \ref{genbtheorem_0}.  
We obtain an \textit{exact}, continuous linear programming reformulation for GB($\gamma$) with 
  $$O(\, 2^{(\omega + 1)L(\gamma)} \, n\, L(\gamma) \, ) \ = \ O( (2 \pi/\epsilon)^{\omega + 1} n \log (\pi/\epsilon) \,)$$
variables and constraints.
In view of Lemma \ref{tobinary}, and the fact that $\delta = 1 - (1 - \gamma)^{\pi} \le \pi \gamma = \epsilon$,  the proof of Theorem \ref{potheorem} is complete.

\subsubsection{Can the dependence on $\epsilon$ be improved upon?}
A reader may wonder why or if ``exact'' 
feasibility (or optimality) for PO cannot be guaranteed.  From a trivial
perspective, we point out that there exist
simple instances of PO (in fact convex, quadratically constrained problems)
where all feasible solutions have irrational coordinates.  Should that be
the case, if any algorithm outputs an explicit numerical solution in finite
time, such a solution will be infeasible.  A different perspective is
that discussed in Example \ref{minknapsackex}.  As shown there
we cannot expect to obtain 
an exact optimal solution in polynomial time, even in the bounded tree-width
case, and even if there is a
rational optimal solution, unless P = NP.

To address either issue one can, instead, attempt to 
output solutions that are approximately feasible.  The approximation
scheme given by Theorem \ref{potheorem} has two characteristics: first, it
allows a violation of each constraint by $\epsilon$ times the 1-norm of the constraint, and second, the running time is pseudopolynomial in $\epsilon^{-1}$. 
One may wonder if either characteristic can be improved.  For example, one
might ask for constraint violations that are at most $\epsilon$, independent
of the 1-norm of the constraints.  However this is not possible even for
a \textit{fixed} value of $\epsilon$, unless P=NP.  For
completeness, we include a detailed analysis of this fact in Section \ref{app:epsilon} of the Appendix. Intuitively, if we were
allowed to approximately satisfy every constraint with an error that does not depend on the data, we could appropriately
scale constraint coefficients so as to obtain exact solutions to NP-hard problems.

Similarly, it is not possible to reduce the pseudopolynomial dependency on $\epsilon^{-1}$ in general. The precise statement is given in  Section \ref{app:epsilon} of the Appendix as well, and the intuitive reasoning is similar: if there was a formulation of size polynomially dependent on $\log(\epsilon^{-1})$
(and not on $\epsilon^{-1}$) we could again solve NP-hard problems in
polynomial time.

\section{Network mixed-integer polynomial optimization problems}\label{graphminpoly}

Here we return to the \textit{network} polynomial optimization problems presented in the introduction, and
provide a proof of Theorem \ref{npotheorem_0}. 
We will first
motivate the technical approach to be used in this proof.  Consider an NPO instance with graph $G$.  For each $u \in V(G)$, $X_u$ denotes the set of variables associated
with $u$ and $\deg_G(u)$ is the degree of $u$.   At each $u \in V(G)$ we have a set of polynomial constraints of the general form
\begin{eqnarray}
\hspace{-.5in}  && \sum_{\{u,v\} \in \delta(u)} p^k_{u,v}(X_{u} \cup X_v) \ \geq \ 0, \quad \quad k = 1, \ldots, K_u\label{typical2}
\end{eqnarray}
associated with $u$, where each $p^k_{u,v}$ is a polynomial  (possibly $K_u = 0$). Without loss of generality we assume that, for every $u$ and $k=1,\ldots, K_u$ no two polynomials $p^k_{u,v}$ $\{u,v\} \in \delta(u)$ in \eqref{typical2} have a common monomial. If that were not the case, we could always combine the common monomials and assign them to a single $p^k_{u,v}$. This allows us to have 
\begin{equation} \label{pk-norm}
\left\|  \sum_{\{u,v\} \in \delta(u)} p^k_{u,v} \right\|_1 =   \sum_{\{u,v\} \in \delta(u)}\left\|   p^k_{ u,v} \right\|_1, \quad k=1,\ldots, K_u.
\end{equation}

Now define
\begin{subequations}\label{Kdef}
\begin{eqnarray}
  \Delta & \doteq & \max_{v \in V(G)}\{ \, |X_v| + K_v \, \}, \mbox{ and} \\
  D & \doteq & \max_{v \in V(G)}\{ \deg_G (v) \}.
\end{eqnarray}
\end{subequations}
As discussed above, we cannot reduce Theorem \ref{npotheorem_0}
to Theorem \ref{potheorem_0} because even if the graph $G$
underlying the NPO problem has small tree-width, the same
may not be the case for the intersection graph of the constraints: a constraint
(\ref{typical2}) can yield a clique of size
$|\delta(u)|$ in the intersection graph (or larger, if e.g. some
$|X_v| > 1$).  In other words,
high degree vertices in $G$ result in large tree-width of the
intersection graph.

\begin{wrapfigure}{r}{0.5\textwidth}
  \begin{center}
      \vskip -15pt
    \includegraphics[height=1in]{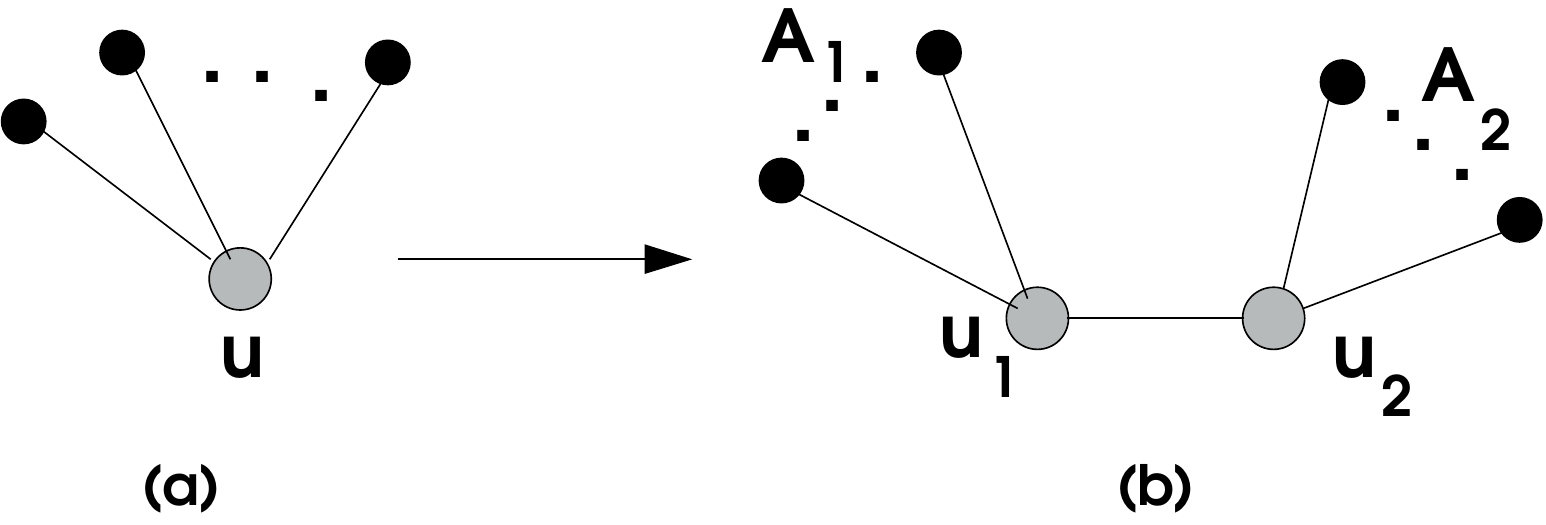}
      \vskip -5pt
    \caption{Vertex splitting.\label{fig:split1}}
  \end{center}
  \vskip -15pt
\end{wrapfigure}

One immediate idea is to employ the technique of \textit{vertex splitting}\footnote{See \cite{G-Splitting} for a column splitting technique used in interior point methods for linear programming.}.  Suppose $u \in V(G)$ has degree larger than three and consider a partition of $\delta(u)$ into two sets $A_1$, $A_2$.
We obtain a new graph $G'$ from $G$ by replacing $u$ with
two new vertices, $u_1$ and $u_2$, introducing the edge
$\{u_1, u_2\}$ and replacing each edge $\{u, v\} \in A_S$ (for $S = 1,2$)
with $\{u_S, v\}$. See Figure \ref{fig:split1}.

\begin{figure}[h] 
\centering
      \includegraphics[height=1.5in]{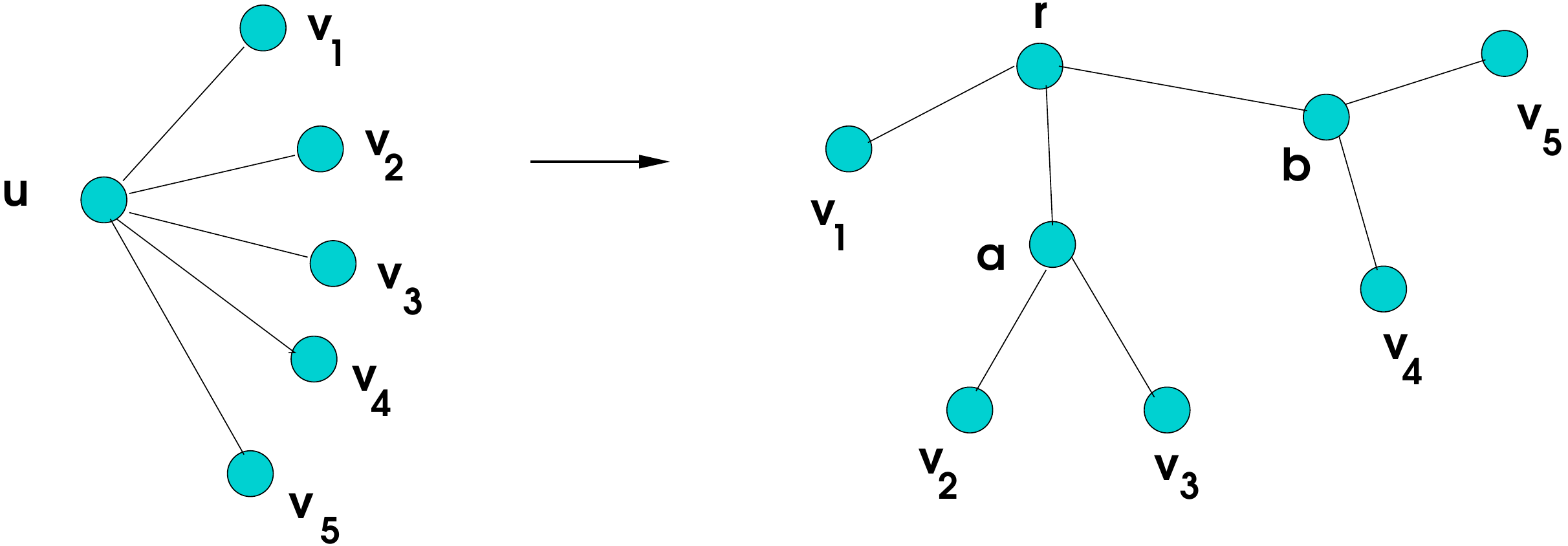}
      \caption{Example of complete vertex splitting.\label{fig:split5}}
\end{figure}

Repeating this procedure, given $u \in V(G)$ with $\deg_G(u) > 3$ we can
replace $u$ and the set of edges $\delta(u)$ with a tree, where each internal vertex will have degree equal to three. To illustrate this construction, consider the example on Figure \ref{fig:split5}. 

\noindent In this figure, a degree-5 node $u$ and the set of edges $\{ (u,v_i) \, : \, 1 \le i \le 5\}$ is converted into a tree with three internal vertices ($r, a$ and $b$) and with each edge $(u,v_i)$ having a corresponding edge in the tree.
To keep the illustration simple, suppose there is a unique constraint of the form 
\begin{eqnarray}
\hspace{-.5in}  && \sum_{i=1}^5 p_{u,v_i}(X_{u} \cup X_{v_i}) \ \geq \ 0. \label{typical5}
\end{eqnarray}
associated with $u$ in the given NPO. This constraint can be transformed into the following system of constraints:
\begin{subequations}\label{split5sys}
  \begin{eqnarray}
        p_{ u,v_2}(X_{u} \cup X_{v_2}) & = & \|p_{u,v_2}\|_1 \, ( y^+_{v_2} \, - \,  y^-_{ v_2})  \label{split5:v2} \\
       p_{ u,v_3}(X_{u} \cup X_{v_3}) & = & \|p_{u,v_3}\|_1 \, ( y^+_{v_3} \, - \,  y^-_{ v_3}) \label{split5:v3} \\
     (\|p_{u,v_2}\|_1  + \|p_{u,v_3}\|_1) (y^+_a - y^-_a) & = & \|p_{u,v_2}\|_1 \, ( y^+_{v_2} \, - \,  y^-_{ v_2}) \ + \ \|p_{u,v_3}\|_1( y^+_{v_3} \, - \,  y^-_{ v_3}) \label{split5:va} \\
    && \nonumber \\
    p_{ u,v_4}(X_{u} \cup X_{v_4}) & = & \|p_{u,v_4}\|_1 \, ( y^+_{v_4} \, - \,  y^-_{ v_4}) \label{split5:v4} \\
    p_{ u,v_5}(X_{u} \cup X_{v_5}) & = & \|p_{u,v_5}\|_1 \, ( y^+_{v_5} \, - \,  y^-_{ v_5}) \label{split5:v5} \\
    (\|p_{u,v_4}\|_1  + \|p_{u,v_5}\|_1) (y^+_b - y^-_b) & = & \|p_{u,v_4}\|_1 \, ( y^+_{v_4} \, - \,  y^-_{ v_4}) \ + \ \|p_{u,v_5}\|_1( y^+_{v_5} \, - \,  y^-_{ v_5}) \label{split5:vb} \\
    && \nonumber \\
    p_{ u,v_1}(X_{u} \cup X_{v_1}) & = & \|p_{u,v_1}\|_1 \, ( y^+_{v_1} \, - \,  y^-_{ v_1}) \label{split5:v1} \\
    && \nonumber \\
\left( \sum_{j = 1}^5 \| p_{u,v_j}\|_1 \right) y^+_r    & = &  \|p_{u,v_1}\|_1 \, ( y^+_{v_1} \, - \,  y^-_{ v_1}) +  (\|p_{u,v_2}\|_1  + \|p_{u,v_3}\|_1) (y^+_a - y^-_a) \nonumber \\
&& \hspace{.2in} + (\|p_{u,v_4}\|_1  + \|p_{u,v_5}\|_1) (y^+_b - y^-_b) \label{split5:vr} \\
&& \nonumber \\
    &&  0 \le y^+, y^- \le 1. \nonumber
  \end{eqnarray}
\end{subequations}

\vspace{.1in}
Clearly, substituting \eqref{typical5} associated with $u$ with the equivalent system \eqref{split5sys} yields an equivalent NPO where 
\begin{eqnarray*}
X_{a} &=& X_u \cup \{ y^+_{a}, y^-_{a}, y^+_{v_2}, y^-_{v_2}, y^+_{v_3}, y^-_{v_3} \}  \\
X_{b} &=& X_u \cup \{ y^+_{b}, y^-_{b}, y^+_{v_4}, y^-_{v_4}, y^+_{v_4}, y^-_{v_4} \} \\
X_{r} &=& X_u \cup \{ y^+_{r}, y^+_{v_1}, y^-_{v_1}\} 
\end{eqnarray*}
 and constraints \eqref{split5:v2}, \eqref{split5:v3} and \eqref{split5:va} are associated with vertex $a$, constraints \eqref{split5:v4}, \eqref{split5:v5} and \eqref{split5:vb} are associated with vertex $b$, and constraints \eqref{split5:v1} and \eqref{split5:vr} are associated to vertex $r$. It is important to notice that we do not associate variables $y^+_{v_i}$ and $y^-_{v_i}$ with nodes $v_i$, but rather with internal vertices of the tree. The motivation for this detail is that vertices $v_i$ may have degree greater than three, and would (later) be split as well. Thus, adding new variables to $X_{v_i}$ might create difficulties when defining the general procedure. This is avoided by keeping the same sets $X_{v_i}$ associated with $v_i$ after a neighbor of $v_i$ is split, as in the previous example.

We indicate the general formal procedure next;
however we warn the reader in advance that this strategy may not directly
deliver Theorem \ref{npotheorem_0} because the splitting process may (if not
chosen with care) produce a graph with much higher tree-width than $G$. This
is a technical point that we will address later (Section \ref{slim}).

To describe the general procedure we use the following notation: given a tree $T$, an edge of $T$ is called \textit{pendant} if it is incident with a leaf, and non-leaf vertex is called \textit{internal}; the set of internal vertices is
denoted by $int(T)$.  Now fix a vertex $u$ of $G$ such that $\deg_G(u) > 3$. Let $\hat T_u$ be an arbitrary tree where
\begin{itemize}
\item $int(\hat T_u) \cap V(G) = \emptyset$, $\hat T_u$ has $\deg_G(u)$ leaves, and for each edge $\{u,v\} \in E(G)$ there is one pendant edge $\{i, v\}$ of $\hat T_u$, and
\item each internal vertex of $\hat T_u$ has degree equal to three,
\end{itemize}
Then \textit{completely splitting $u$ using $\hat T_u$} yields a new graph, $\breve G$ where $V(\breve G) \ = \ ( V(G) - u ) \cup int(\hat T_u)$, and $E(\breve G) \ = \ ( E(G) - \delta_G(u) ) \cup E(\hat T_u)$.

To obtain an NPO system in $\breve G$ equivalent to the original
system, we replace each constraint (\ref{typical2}) that is associated with
$u$ with a family of constraints associated with the internal vertices of
$\hat T_u$. To do so, pick an arbitrary non-leaf vertex $r$ of $\hat T_u$ and view
$\hat T_u$ as \textit{rooted} at $r$ (i.e., oriented away from $r$), and define, for
each internal vertex $i$ of $\hat T_u$ and $k=1,\ldots, K_u$,
\begin{eqnarray}
  \nu_i^k & = & \sum \left\{ \| p^k_{u,v}\|_1 \, : \, (u,v) \in \delta(u), \, \mbox{$v$ a descendant of $i$ in $\hat T_u$} \right\} \nonumber
\end{eqnarray}    
Then, for each internal vertex
$i$ of  $\hat T_u$ we have
\begin{itemize}
\item [(a)] All variables in $X_u$ are associated with $i$.
\item [(b)] For each $k=1,\ldots, K_u$,
 if $i \neq r$ we additionally associate two variables, $y^+_{i,k}, ~ y^-_{i,k} \in [0,1]$,  and if $i = r$ we only add $y^+_{r,k} \in [0,1]$ with $i$.  If $i$ has a child $v$ that is a \textit{leaf}, and hence
  $\{u,v\} \in \delta(u)$, then we associate two
  additional variables $y^+_{v,k}, ~ y^-_{v,k} \in [0,1]$ with $i$.
\item [(c)] If $i \neq r$, then letting $j, l$ be its children we write the following constraints associated with $i$:
  \begin{eqnarray}
    &&    \nu^k_j (y^+_{j,k} - y^-_{j,k}) + \nu^k_l (y^+_{l,k} - y^-_{l,k}) \, = \, \nu^k_i (y^+_{i,k} - y^-_{i,k}) \quad \quad 1 \le k \le K_u. \label{internal}
    \end{eqnarray}
      If on the other hand $i = r$, let $j, l, h$ be its children. Then
      we write
  \begin{eqnarray}
    &&    \nu^k_r y^+_{r,k} - \sum_{s = j, l, h} \nu^k_s (y^+_{s, k} - y^-_{s, k})\  = \ 0  \quad \quad 1 \le k \le K_u. \label{root}
  \end{eqnarray}
\end{itemize}
Finally, given a leaf $v$ of $\hat T_u$ with parent $i$, then by construction
$\{u,v\} \in \delta(u)$. Then we add the following constraints,
associated with $i$ (not $v$), 
  \begin{eqnarray}
    &&    \|p^k_{u,v}\|_1 (y^+_{v,k} - y^-_{v,k}) \ - \ p^k_{u,v}(X_u \cup X_v) \ = \ 0 \quad \quad 1 \le k \le K_u.   \label{leaf}
  \end{eqnarray}

%

Let us denote the initial NPO by $P$.  Since the sum of constraints (\ref{internal}), (\ref{root}) and (\ref{leaf}) is
(\ref{typical2}) we have obtained an NPO equivalent to $P$.  Note that for 
any $v \in V(G) - u$ the degree of $v$ is unchanged, as is the set $X_v$ and the
set of constraints associated with $v$. Thus,
proceeding in the above manner
with every vertex $u \in V(G)$ with $\deg_G(u) > 3$ we will obtain the final
graph (which we denote by $G'$) of maximum degree $\le 3$ and an NPO, denoted by $P'$ which is equivalent
to $P$. 

The following lemmas lay out the strategy that we will follow
to prove Theorem \ref{npotheorem_0}. To prepare for these, we need a technical remark, which follows from
the definition of tree-decomposition.
\begin{remark} \label{connected} Suppose that $(T,Q)$ is a tree-decomposition of a graph $H$.
  Suppose a subset of vertices $S$ induces a connected subgraph of $H$.  Then
  $$ \{t \in V(T) \, : \, Q_t \cap S \neq \emptyset\}$$
  induces a subtree of $T$.\end{remark}

Returning to our construction,  
we assume that we have a tree decomposition $(T',Q')$ of $G'$.  Note that
each vertex $v$ of $G'$ is derived from some vertex, say $u(v) \in V(G)$ and the set of variables associated with $v$ under $P'$ is a copy of
$X_{u(v)}$, together with the $y^+$ and $y^-$ variables introduced above.  We have
one pair of such variables per each constraint (\ref{typical2}) associated with
$u(v)$ and an extra pair when $v$ is the parent of a leaf in $\hat T_{u(v)}$.

We consider now the following pair: $(T', \tilde Q')$, with each $\tilde Q'_t$ $t\in T'$ defined as follows. For each $v \in Q'_t$, the set $\tilde Q'_t$
will include (1) $X_{u(v)}$ together with the $y^+$ and $y^-$ variables associated
with $v$ or its children in $\hat T_{u(v)}$. This can include up to $|X_v| + 7K_v$ variables for each $v$ (the bound is tight at the root $r$ which is the only internal vertex with three children; further we are counting $y^+_r$). (2) If $v$ is the parent of leaf in $\hat T_{u(v)}$, say $w$, then we also add to $\tilde Q'_t$ the set $X_{u(w)}$. Further, the number of constraints associated with
internal vertices is clearly $\le 3K_v$.

\begin{lemma} The pair $(T', \tilde Q')$ is a tree-decomposition of the intersection graph of $P'$.\end{lemma}
\noindent {\em Proof.} We need to show that (a) for any variable $x_j$ of $P'$.  
the set of vertices $t$ of $T'$ such that $\tilde Q'_t$ contains $x_j$ forms a subtree of $T'$, and (b) that for any edge $\{x_k, x_j\}$ of the intersection graph
of $P'$ the pair $x_k, x_j$ are found in a common set $\tilde Q'_t$.

Part (b) follows directly from the construction of the sets $\tilde Q'_t$, as each of the constraints \eqref{internal}, \eqref{root} and \eqref{leaf} involve only one internal node and its children, including the case when some children are leaves, which is accounted for in (2).

As for part (a), the statement
is clear if $x_j$ is one of the $y^+, y^-$ variables. If, instead,
$x_j$ is contained in some set $X_{u(v)}$ then $x_j$ is associated with every vertex in $\hat T_{u(v)}$.  By definition of NPOs, the set of $w$ vertices of $G$ such
that $x_j \in X_w$ forms a connected subgraph of $G$.  It follows that the set
of vertices $v$ of $G'$ such that $x_j$ is associated with $v$ in $P'$ forms a
connected subgraph of $G'$.  Part (a) now follows from Remark \ref{connected}. \QED

\begin{lemma} Suppose that $G'$ has tree-width $ \le W$.  Then
  the intersection graph of $P'$ has tree-width $\le 7 \Delta (W + 1) - 1 $. \end{lemma}
\noindent {\em Proof.} Consider a tree-decomposition $(T',Q')$ of $G'$ of
width $\le W$ and construct $\tilde{Q}'$ as before.  We claim that the width of
$(T', \tilde Q')$ is at most $7 \Delta (W + 1) - 1$. But this is clear since $|Q'_t| \leq W+1$, and each $v\in Q'_t$ will contribute at most three extra sets of the type $X_w$ for some $w$, along with the $|X_v| + 7K_v$ quantities stated before. Hence,
$$|\tilde Q'_t| \leq 7 \Delta (W + 1)$$ 
as desired.
\QED

\begin{lemma} \label{hehehe}
  Given $0 < \epsilon < 1$, there is a linear programming formulation LP for the NPO problem on $G$, of size $O((D \pi/\epsilon)^{O(\Delta W)}\, n \, \log \epsilon^{-1})$, that solves the problem within scaled tolerance $\epsilon$.
\end{lemma}
\noindent {\em Proof.} Suppose we apply Theorem \ref{potheorem} to problem
$P'$, for a given tolerance $0 < \theta < 1$. Note that, in order to do so rigorously, equations \eqref{internal} and \eqref{leaf} must be first transformed to 2 inequality constraints each. 

The resulting linear programming formulation, LP, will have size $O((2\pi/\theta)^{O(\Delta W)}\, n \, \log \theta^{-1})$ and will yield $\theta$-approximate solutions to $P'$.  In particular,
these solutions may violate a constraint of $P'$ by an amount proportional
to the $1$-norm of the constraint, times $\theta$, as guaranteed by
Theorem \ref{potheorem}\footnote{In the case of a constraint that had the form $p(x) = 0$ in $P'$, and which was expressed as two inequalities, it means that a $\theta$-scaled feasible vector $\hat{x}$ will satisfy
$$p(\hat{x}) \geq -\theta \|p\|_1 \ \ \mbox{and} \  - p(\hat{x}) \geq -\theta \|p\|_1 \ \mbox{and thus} \ |p(\hat{x})| \leq \theta \| p\|_1.$$}.

Consider a constraint (\ref{typical2}) of problem $P$, associated with a vertex $u$ of $G$, with 1-norm denoted by $\|p_u^k\|_1$.
In $P'$ this constraint has been replaced by a family $\F$ of constraints, one
per each vertex of $\hat T_u$ whose sum yields (\ref{typical2}).  Simple algebra shows that the solution to LP may violate (\ref{typical2}) by an amount equal to the sum of violations of the constraints in $\F$.  The $1$-norm of each such constraint is
at most  $4 \|p_u^k\|_1$, since if $i$ is an internal node of some tree $\hat T_u$ with children $j,l$ then 
$$\nu_i^k = \nu_j^k + \nu_l^k $$ 
and for a root node $r$ with children $j,h,l$ 
$$\|p_u^k\|_1 = \nu_r^k = \nu_j^k + \nu_h^k +\nu_l^k $$
by \eqref{pk-norm}. Since $\hat T_u$ has $\deg_G(u)$ leaves and all internal vertices have degree 3, it can be easily shown that it must have less than $2 \deg_G(u)$ vertices in total. From this it follows that the LP solution violates (\ref{typical2}) by
at most $8 \|p_u^k\|_1 \deg_G(u) \theta$.  Choosing $\theta = (8D)^{-1} \epsilon $, and noting that for a constant $c$
$$ (c D\pi/ \epsilon)^{O(\Delta E)} \in O( D\pi/ \epsilon)^{O(\Delta E)})$$
yields the desired result. \QED

\vspace{.1in}

As per Lemma \ref{hehehe}, these constructions obtain Theorem \ref{npotheorem_0} provided
that we find a  vertex splitting with  $W = O(\omega)$. As discussed
in the next section, this condition may fail to hold under an arbitrary vertex
splitting, and some care is needed.

\subsubsection{Finding good vertex splittings}\label{slim}
\begin{figure}[h] 
\centering
\includegraphics[height=1.5in]{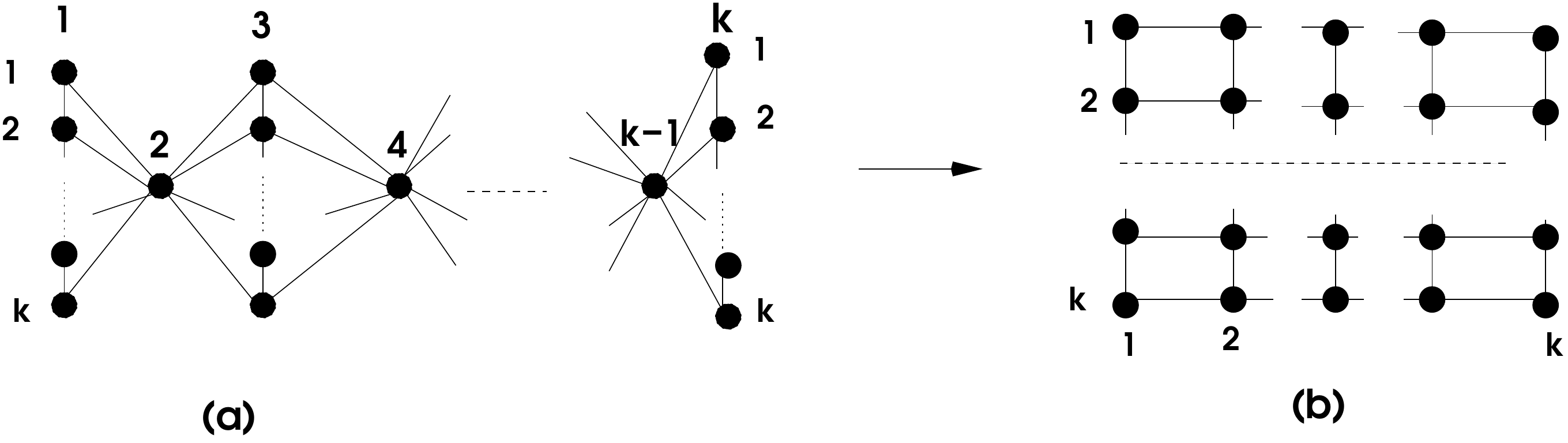}
\caption{Incorrect vertex splitting.\label{fig:split2}}
\end{figure}

Consider the example given in Figure
\ref{fig:split2}.
Here the vertices of the graph shown in (a) are arranged into $k > 1$ columns.
The odd-numbered columns have $k$ vertices each, which induce a path, while the
even-numbered columns have a single vertex which is adjacent to
all vertices in the preceding and following columns. It can be shown
that this graph has tree-width $3$.  In (b) we show the outcome
after splitting vertices so that the maximum degree is four.
This second graph has tree-width $k$, and further splitting the
degree-four vertices will not change this fact.

\begin{figure}[h] 
\centering
\includegraphics[height=1.1in]{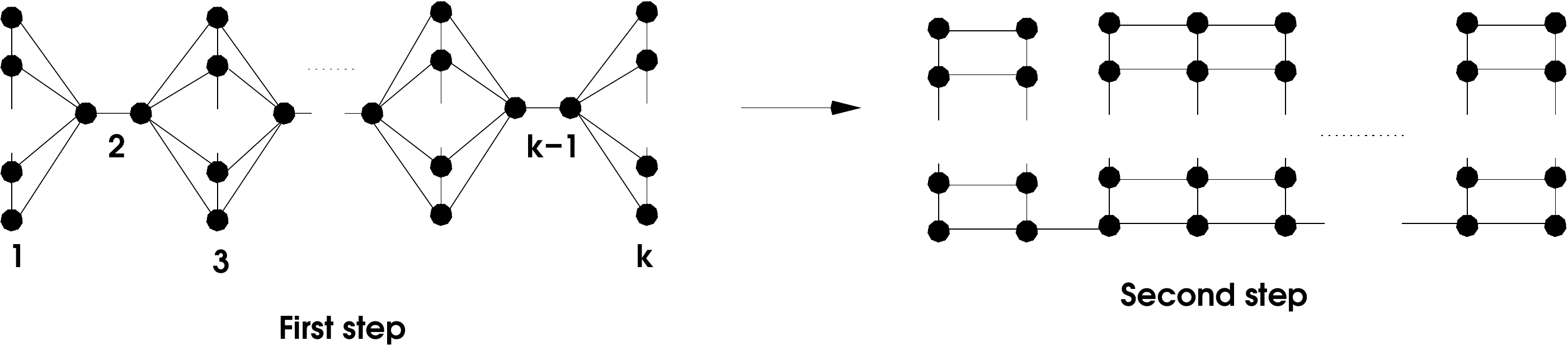}
\caption{Correct vertex splitting of graph in Figure \ref{fig:split2}(a).\label{fig:split3}}
\end{figure}

In contrast to this situation, suppose that we
split the graph in 
Figure \ref{fig:split2}(a) in two steps as shown in
Figure \ref{fig:split3}.  The tree-width of the final graph is also 3.  The difference between Figures \ref{fig:split2} and \ref{fig:split3} is explained
by the fact that the splitting initiated by the ``first step'' in Figure \ref{fig:split3} exploits the
tree-decomposition of width $3$ of the graph Figure \ref{fig:split2}(a).

Next we turn to a formal approach that produces the desired outcome in the general
setting.  Given a graph $G$, a \textit{simplification} of $G$ will be graph $\bar G$ obtained by a sequence of complete vertex splittings, such that the maximum degree of
a vertex in $\bar G$ is $\le 3$. The following Lemma will show how to obtain a simplification of a graph via a vertex splitting that maintains tree-width up to a constant factor. In the proof, the trees $\hat T_u$ that yield the general splitting procedure stated above, will be explicitly defined.

\begin{lemma} \label{structural} Let $G$ be an undirected graph 
and $(T, Q)$ a tree-decomposition of $G$ of width $Z$.  Then there is
a simplification $\bar G$ of $G$ and a tree-decomposition
$(\bar T, \bar Q)$ of $\bar G$ of width at most $2 Z+1$.
\end{lemma}
\noindent {\em Proof.} We first modify $(T, Q)$ in a sequence of steps.\\

\noindent {\bf Step 1.} For any edge $e = \{u, v\} \in E(G)$, choose an arbitrary 
$t \in V(T)$ with $e \subseteq Q_t$.  Then we modify $T$ by adding to $T$ a new vertex, $t^e$ and the edge $\{t^e, t\}$.  Further, we set $Q_{t^e} = \{u,v\}$.  \\

\noindent {\bf Step 2.} Without loss of generality, every vertex of $T$ has degree
at most $3$.  To attain this condition, consider any $t \in V(T)$ with
$\delta_T(t) = \{s_1, \ldots, s_d\}$ (say) where $d > 3$.  Then we alter $T$ by
replacing $t$ with two vertices adjacent vertices $t^1$ and $t^2$, such that $t^1$ is also adjacent to $s_1$ and $s_2$ and $t^2$ is adjacent to $s_3, \ldots, s_d$.  
Finally, we set $Q_{t^1} = Q_{t^2} = Q_t$. Continuing inductively we will attain
the desired condition.\\

\noindent {\bf Step 3.} For any vertex $u \in V(G)$ let $T_u$ be 
the subtree of $T$ consisting of vertices $t$ with $u \in Q_t$, and
$\breve T_u$ be the subtree of $T_u$ that spans $\{t^e \, : \, e \in \delta_G(u)\}$ (which is a subset of the leaves of $T_u$).  Then we modify $(T,Q)$ by replacing $T_u$ with $\breve T_u$, yielding
a new tree-decomposition of same or smaller width.  In other words, without loss of generality \textit{every} leaf of $T_u$ is of the form $t^e$ for some $e \in \delta_G(u)$.\\

We can now describe our vertex splitting scheme.  Consider 
$u \in V(G)$ with $\deg_G(u) > 3$.  We say that a vertex of $T_u$ is \textit{blue} if it is either a leaf or of degree three in $T_u$.  Now we form the tree
$\hat T_u$ whose vertex-set is the set of blue vertices of $T_u$, and whose edge-set is obtained as follows.  By construction, $E(T_u)$ can be partitioned into a set of
paths whose endpoints are blue and which contain no other blue vertices.  For
each such path, with endpoints $a$ and $b$ (say), the tree $\hat T_u$ will contain the edge $\{a,b\}$ (in other words, $T_u$ can be obtained from $\hat T_u$ by
subdividing some edges and so $T_u$ and $\hat T_u$ are topologically equivalent). Note that $\hat T_u$ has $\deg_G(u)$ leaves, each internal vertex with degree 3, and for each edge $\{u,v\} \in E(G)$ there is one pendant edge, as needed.

Let $\breve G$ be the graph obtained by the complete splitting of $u$ using
$\hat T_u$. For each internal vertex $t\in V(\hat T_u)$ we name $u^t$ the corresponding new vertex in $\breve G$, to emphasize that each non-leaf vertex in $\hat T_u$ will create a copy of $u$ (recall that the leaves in $\hat T_u$ will correspond to the neighbors of $u$). This operation does not change the degree of any vertex
$v \in V(G)$ with $v \neq u$.  The eventual graph $\bar G$ in the proof will be obtained by applying complete splittings at every vertex of degree $> 3$ in $G$.

Returning to $\breve G$ we construct a tree decomposition $(T, \breve Q)$ as follows.  First, let us regard the tree $T_u$ as rooted at some internal blue vertex $r(u)$.  For a vertex $t \in V(T_u)$ let $R_u(t)$ be the closest blue \textit{ancestor} of $t$ in $T_u$; we write $R_u(r(u)) = r(u)$.  Then, for $t \in V(T)$, we set
\begin{eqnarray} \label{splitit}\breve Q_t \ = \  \left\{ \begin{array} {ll}
		(Q_t - u) \cup \{u^{R_u(t)}\}, & \mbox{if $t \in V(T_u)$ and $\deg_{T_u}(t) = 2$ }, \\
    (Q_t - u) \cup \{u^t, u^{R_u(t)}\}, & \mbox{if $t \in V(T_u)$ and $\deg_{T_u}(t) \neq 2$}, \\
    Q_t, & \mbox{if $t \notin V(T_u)$.} 
		\end{array} 
	\right.
\end{eqnarray}
Now we argue that $(T, \breve Q)$ is a tree-decomposition of $\breve G$. To see this, note that if $t \in V(\hat T_u)$ then $u^t \in \breve Q_s$ iff $s = t$ or $s$ is a child of $t$ in $\hat T_u$, thus the endpoints of any edge $\{u^t, u^s\}$, where $t$ is the parent of $s$, will be contained in $\breve Q_s$.
Further, for any edge of $e = \{u, v\}$ of $G$, by Step 3 above there
will be a leaf $t^e$ of $T_u$ such that the edge $\{t^e, R_u(t^e)\} \in E(\hat T_u)$. This corresponds to a pendant edge $\{u^{R(t^e)}, v\} \in E(\breve{G})$ and by construction both $v \in \breve Q_{t^e}$ and $u^{R(t^e)} \in \breve Q_{t^e}$. The fact that every vertex in $\breve{G}$ induces a connected subgraph in $T$ can be easily verified. This completes the argument $(T, \breve Q)$ is a tree-decomposition of $\breve G$.

Notice that for $v \in V(G)$ with $v \neq u$, the subtree $T_v$ is the same
in $(T, Q)$ and $(T, \breve Q)$.  Thus, applying the complete splitting of
every vertex of $G$ of degree greater than three, and modifying the tree-decompostion as in (\ref{splitit}) will produce a tree-decomposition $(T, \bar Q)$ of
the final graph $\bar G$.

By construction, for each $t \in V(T)$ we obtain $\bar Q_t$ from $Q_t$ by replacing each element with (at most) two new elements.  Thus, since $|Q_t| \leq Z + 1$, the width of $(T, \breve Q)$ is at most $2(Z+1) - 1$. \QED \\


\section{Pure binary problems}\label{purebinary}
In this section we will consider Theorem \ref{genbtheorem_0} of the Introduction.
As we mentioned above, it is one of the building 
blocks towards the other main results, but is of independent interest as well. We will provide
additional background, a deep analysis of this result, and state and prove an expanded
version of the Theorem. First we begin with some examples
for problem GB.

\begin{example} \label{ex0} (Linear binary integer programming).  Let $A$ be an
$m \times n$ matrix, and consider a problem $\min \{ c^T x \, : \, Ax \ge b, \ x \in \{0,1\}^n\}$.  To view this problem as a special case of GB, we set for $1 \le i \le m$,
$K[i] = \{ 1 \le j \le n \, : \, a_{ij} \neq 0 \}$ and $S[i] = \{ x \in \{0,1\}^{K[i]} \, : \, \sum_{j \in K[i]} a_{ij} x_j \ge b_i \}$. 
\end{example}
 In this special case, problem GB
can be addressed by a variety of methods.  Of particular interest in this paper are
the reformulation or lifting methods of \cite{Lovasz91conesof} and \cite{sheraliadams}.
Next we consider a more complex example, chosen to highlight
the general nature of the problem.
\begin{example} \label{firstexample} Let $d, r, p$ be positive integers. Consider a
constrained semidefinite program over binary variables of the form
\begin{subequations} \label{sdpbin}
\begin{eqnarray}
&&  \min ~ \sum_{k = 1}^r \sum_{i = 1}^{d} \sum_{j = 1}^{d}c_{kij} X^k_{i,j} \\
\mbox{subject to:} && \ M^k \bullet X^k \ = \ b_k, \quad 1 \le k \le r, \label{bullet}\\
&& \ X^k \ \in \ S^+_{d}, \quad 1 \le k \le r, \label{semi}\\
&& \ \sum_{i,j} X^k_{i,j} \ \equiv  0 \hspace{-5pt} \mod p, \quad \ 1 \le k \le r, \label{modulo}\\
&& \ X^k_{i,1} \ = \ X^{k-1}_{i,d}, \quad 1 \le i \le d, \quad 2 \le k \le r, \label{overlap}\\
&& \ X^k_{i,j} \, \in \, \{0,1\}, \ \forall i, j, k.
\end{eqnarray}
\end{subequations}

\noindent Here $S^+_{d}$ is the set of $d \times d$ positive-semidefinite matrices, $M_1, \ldots, M_r$ are
symmetric $d \times d$ matrices, and $b$ and $c$ are vectors.  Constraint 
(\ref{overlap}) states that the first column of matrix $X^k$ is identical
to the last column of matrix $X^{k-1}$.  

We obtain an instance of problem GB
with $m = 2r - 1$, as follows.  First, for each $1 \le k \le r$ we let $K[k]$ be the
set of triples $(i,j,k)$ with $1 \le i, j \le r$, and $S^k$ to be the set of 
binary values
$X^k_{i,j}$ that satisfy (\ref{bullet})-(\ref{modulo}).  Next, for each 
$2 \le k \le r$ we let $K[r + k -1]$ be the set of all triples $(i,1,k-1)$ and all triples $(i,d,k)$ and $S^{r + k -1}$ to be the set of binary values (indexed by $K[r + k -1]$) such
that (\ref{overlap}) holds.
\end{example}
In the case
of this example, a direct application of standard integer programming methods
appears difficult. Moreover, we stress that the sets $S^i$ in problem GB are completely generic
and that the membership oracle perspective can prove useful as we discuss below.

Theorem \ref{genbtheorem_0} concerns the tree-width of the intersection graph of a
problem of type GB.  Recall that as 
per Definition \ref{intersectgraph}, given a problem instance $\I$ of GB, 
the intersection graph for $\I$ has a vertex for each $1\leq j \leq n$, and an
edge $\{j, k\}$ whenever there exists $1 \le i \le m$ such that 
$ \{j, k\} \ \subseteq \  K[i]$, that is to say, $j$ and $k$ appear in a common constraint in problem GB.

\begin{example}
(Example \ref{firstexample}, continued). Here the set of variables is given by
$$ \{ (i,j,k) \, : \, 1 \le k \le r \, \mbox{ and } \  1 \le i, j \le d\}. $$
The intersection graph of the problem will have 
\begin{itemize} 
\item [(a)] the edge $\{(i,j,k), (i',j',k)\}$ for all $1 \le k \le r$ and $1 \le i, j, i', j' \le d$, arising from constraints (\ref{bullet})-(\ref{modulo})
\item [(b)] the edge $\{(i,1,k), (i,d,k-1)\}$ for each $1 \le k < r$ and $1 \le i \le d$,
arising from constraints (\ref{overlap}).
\end{itemize}
A tree-decomposition $(T, Q)$ of the intersection graph, of width $O(d^2)$, is 
obtained as follows.  Here,
$T$ is a path with vertices $v_1,u_2,v_2,u_3,  \ldots, v_{r-1}, u_r, v_r$.  For $1 \le k \le r$ we set $Q_{v_k} = \{(i,j,k) \, : \, 1 \le i, j \le d\}$ and for $2 \le k \le r$ 
we set $Q_{u_k} = \{ (i,1,k), (i,d,k-1) \, : \, 1 \le i \le d \}$.  Sets $Q_{v_k}$ account
for all edges of type (a), whereas the sets $Q_{u_k}$ cover all edges of type (b).  Thus Theorem \ref{pbtheorem} states that there
is an LP formulation for problem (\ref{sdpbin}) with $O(2^{d^2} d^2 r)$ variables and
constraints.
\end{example}

We now state the main result we will prove 
regarding problem GB,  which implies Theorem \ref{genbtheorem_0} (a).  A proof
of part (b) of Theorem \ref{genbtheorem_0} is given in Section \ref{counting} of the Appendix.

\begin{theorem} \label{pbtheorem} Let 
$(T, Q)$ be a tree-decomposition of
the intersection graph of a problem GB.  Then 
\begin{itemize}
\item[(a)] There is an exact (continuous) linear programming reformulation 
 with $O(\sum_t 2^{|Q_t|})$ variables and constraints, the
same objective vector $c$ and constraints with $\{0,1,-1\}$-valued coefficients.
\item[(b)] The formulation can be constructed by performing $2^{\omega}m$ oracle queries and with 
additional workload $O^*(\omega 2^{\omega} \sum_{t \in V(T)} |\left\{i \, : \, K[i] \subseteq Q_t \right\}| + \omega mn)$, where the ``*'' notation indicates logarithmic factors in $m$ or $n$.
\end{itemize}
As a corollary, if the width of $(T, Q)$ is $\omega$, the formulation 
has $O(2^{\omega} n)$ variables and constraints.  Hence for each fixed $\omega$ the formulation has linear size.
\end{theorem}
The ``corollary'' statement follows because if an $n$-vertex graph has a tree-decomposition
of width $\omega$, say, then it has one with the same width and where in addition the tree has at most $n$ vertices (see Remark \ref{tdecompsize}, below).  To illustrate, we show what this result implies when applied to one of our previous examples:

\subsubsection{Remark. Reduction to the linear case}\label{alllinear}
Consider a problem instance of GB. An apparently simpler 
alternative to the general approach we follow would be to construct, for
$1 \le i \le m$, the polyhedron 
$$ P_i \ \doteq \conv\left\{ x \in \{0,1\}^{K[i]} \, : \, x  \in S^i \right\} \ \subseteq \ \R^{K[i]}.$$
Thus we can write $P_i$ as the projection onto $\R^{K[i]}$ of a polyhedron $\{ x \in [0,1]^n \, : \, A^i x \, \ge \, b^i \}$
where each row of $A^i$ has zero entries on any column not in $K[i]$.  
Thus the GB problem can be restated as the equivalent linear \emph{integer} program
\begin{subequations} \label{badLP}
\begin{eqnarray}
 && \min ~ c^Tx \\
\mbox{subject to:} && \quad \quad A^i x \, \ge \, b^i, \quad 1 \le i \le m \label{poly1}\\
&& \quad \quad x \ \in \ \{0,1\}^{n}.
\end{eqnarray}
\end{subequations}
Switching to this formulation makes it possible to apply general integer programming
methods to problem GB. However, this analysis ignores the size of formulation (\ref{badLP}). In particular, for integer $d \ge 1$ large enough there exist examples of $0/1$-polytopes in $\R^d$ with at least $$\left(\frac{d}{\log d}\right)^{d/4}$$ facets (up to constants). See \cite{BaranyPor}, 
\cite{gatzouras05:_lower}, \cite{Z52}.  Using this observation, one can construct
examples of problem GB where the tree-width of the intersection graph is $\omega = d-1$
and each of the matrices $A^i$ has more than
$\omega^{\omega/4}$ rows (see Example \ref{badLPex}, below). This dependence on $\omega$ makes any classical integer programming method more computationally 
expensive than using the method presented above.

\begin{example} \label{badLPex} Choose $d \ge 2$ large enough so that there is a $0/1$-polyhedron $P \subseteq \R^d$ with more than $(c d/\log d)^{d/4}$
facets for some $c$.  Let $P$ be given by the system 
$Ax \ge b$, where $A$ is $M \times d$ ($M \ge (c d/\log d)^{d/4}$).  Choose
$N \ge 1$, and consider the system of inequalities
over binary variables $x^i_j$, for $1 \le i\le N$ and $1 \le j \le d$:
\begin{subequations}
\label{firstpoly}
\begin{eqnarray}
&& A x^i \ge  b, \quad 1 \le i\le N, \label{firstpolyh}\\
&& \quad x^1_j =  x^i_j \quad 2 \le i \le N, \ 1 \le j \le \lfloor d/2 \rfloor. \label{secondpolyh}\\
&& \quad x^i_j \ \mbox{binary for all $i$ and $j$.} 
\end{eqnarray}
\end{subequations}
Constraint (\ref{firstpolyh}) indicates that this system includes $N$ copies
of polyhedron $P$, with each copy described using a different coordinate system.  Constraint
(\ref{secondpolyh}) states  that the first $\lfloor d/2 \rfloor$ coordinates take equal value 
across all such systems.  

Any linear program over (\ref{firstpoly}) is 
can be viewed as an example of problem GB with $m = 2N - 1$; for 
$1 \le i \le N$, $K[i]$ is used to represent the $d$ variables $x^i_j$ ($1 \le j \le d$) and $S^i$ is a copy of the set of binary points contained in $P$ (i.e.
the extreme points of $P$). 

The 
intersection graph of this instance of GB will be the union of $N$ cliques (one for each set of
variables $x^i$) plus the set of edges $\{x^1_1, x^i_1\}$ for $2 \le i \le N$. 
A tree-decomposition $(T, Q)$ of this graph, of width $d-1$, is as follows:  $T$ has vertices
$u(0)$, as well as $u(i)$ and $v(i)$, for $1 \le i \le N$.  Further, $Q_{u(0)} = \left\{ x^1_j \, : \, 1 \le j \le \lfloor d/2 \rfloor \right\}$; and for
$1 \le i \le N$ $Q_{u(i)} = Q_{u(0)} \cup \left\{ x^i_j \, : \, 1 \le j \le \lfloor d/2 \rfloor \right\}$ and $Q_{v(i)} = \{x^i_j, \, 1 \le j \le d\}$.  Thus, $\omega = d$ and Theorem \ref{pbtheorem} states that any linear objective
problem over constraints (\ref{firstpoly}) can be solved as a continuous LP with
$O(2^d d N)$ variables and constraints.  In contrast, 
system (\ref{firstpoly}) has more than $ (c d/\log d)^{d/4} N$ constraints, and in
particular the same is true for the level-$d$ RLT reformulation.
\end{example}
As the example shows, formulation (\ref{badLP}) may 
be exponentially larger than the linear program  stated in Theorem \ref{pbtheorem}.

\subsection {Proof of Theorem \ref{pbtheorem}}  \label{pbproof}
In this section we discuss a construction that yields Theorem \ref{pbtheorem} by
relying on methods from \cite{Laurentima}.  Fundamentally the construction  employs
the ``cone of set functions'' approach of \cite{Lovasz91conesof}
(also see \cite{Laurent01acomparison}), together with
an appropriate version of the ``junction tree theorem'' \cite{graphical} as
developed in \cite{Laurentima}. In addition, we 
provide a second formulation in Section \ref{pbproof2} of the Appendix, and a direct proof.

Consider an instance of problem GB. Let $G$ be the corresponding intersection graph, and $(T, Q)$ be a tree-decomposition of $G$ of width $\omega$. We begin with some
general remarks.

\begin{remark}\label{cliques} Suppose that $(\bar T, \bar Q)$ is a tree-decomposition of a graph $\bar G$.  Then for any clique $K$ of $\bar G$ there exists
$ t \in V(\bar T)$ with $K \subseteq \bar Q_t$. \end{remark}

\noindent As a result, for $1 \le i \le m$ there exists $t \in V(T)$ with $K[i] \subseteq Q_t$, i.e, the indices of the support of each constraint must be contained in some node of the tree-decomposition.

\begin{remark}\label{tdecompsize} Without loss of generality, $|V(T)| \le n$.  To
see this, note that the tree-decomposition $(T,Q)$ gives rise to a chordal
\textit{supergraph} $H$ of $G$.  Since $H$ is chordal, there exists a vertex $u$ whose
neighbors (in $H$) induce a clique.  The claim follows by induction applied to the graph $H - u$, using Remark \ref{cliques} and noting that a tree-decomposition of $H$ is
also a tree-decomposition of $G$.
\end{remark}

\begin{Def}\label{bigdef} Let $t \in V(T)$.
\begin{itemize} 
\item[(a)] We say that $v \in \{0,1\}^{Q_t}$ is \bmath{Q_t}-{\bf feasible} if $v_{K[i]} \in S^i$ for every $1 \le i \le m$ such that $K[i] \subseteq Q_t$.
\item[(b)] Write \bmath{F_t} $ = \, \{ v \in \{0,1\}^{Q_t} \, : \, \mbox{$v$ is $Q_t$-feasible} \}$.
\end{itemize}
\end{Def}




Now we can present the formulation. The variables are as follows:
\begin{itemize}
\item A variable $\lambda^t_v$, for each $t \in V(T)$ and each vector $v \in F_t$.
\item A variable $Z_S$, for each $S \in 2^{V(T)}$ such that $S \subseteq Q_t$ for
some $t \in T$.
\end{itemize}
We also write
$$ V \ = \ V(G) \ = \ \mbox{the set of all variables indices for the given instance of GB},$$
so that $|V| = n$. The formulation is as follows:
\begin{subequations} \label{LPzeta}
\begin{eqnarray}
\mbox{(LPz):} && \min ~ \sum_{j =1}^n c_j Z_{\{j\}} \\
\mbox{s.t.} \quad \forall t \in V(T):&& \nonumber \\
 \quad Z_S & = & \sum\{ \lambda^t_v \, : \, v \in F_t, \ v_j = 1 \ \forall j \in S \} \quad \quad \forall S \subseteq Q_t \label{consistency0}\\
&& \nonumber \\
\sum_{v \in F_t} \lambda^t_v & = & 1, \quad \lambda^t \ge 0. \label{tconvexity0}
\end{eqnarray}
\end{subequations}
Constraints (\ref{consistency0}) enforce consistency across different $t \in V(T)$.  In fact
the $Z$ variables can be eliminated with (\ref{consistency0}) replaced with
relationships among the $\lambda$ variables. Constraint (\ref{consistency0}) can be restated in a more familiar way.  
Given $t \in T$,  (\ref{consistency0}) states:
\begin{eqnarray}
&& Z_S \ = \ \sum_{v \in F(t)} \lambda^t_v \zeta^{supp(v, Q_t)}, \quad \quad \forall S \in 2^{Q_t}. \label{LSt}
\end{eqnarray}
Here, given a set $Y$ and a vector $w \in \R^Y$, $supp(w,Y) = \{j \in Y \, : \, w_j \neq 0 \}$, and for any set $Y$ and $p \subseteq Y$ the vector $\zeta^{p} \in \{0,1\}^{2^Y}$ is defined by setting, for each $q \subseteq Y$,
$$\zeta^{p}_q \ = \ \left\{\begin{tabular}{cl}
1 & \mbox{if $q \subseteq p$}\\
0 & \mbox{otherwise.}
\end{tabular}\right.$$

Constraints (\ref{consistency0})-(\ref{tconvexity0}) describe the Lov\'{a}sz-Schrijver approach
to lifted formulations, restricted to a given set $Q_t$.  It is clear that LPz amounts to a relaxation for the given
problem GB, in the sense that given $\hat x$ feasible for GB then there is a vector $(\hat Z, \hat \lambda)$ feasible for LPz where $\hat Z_{\{j\}} = \hat x_j$ for $1 \le j \le n$. To do so, let $t \in V(T)$ and denote by $\hat x^t$ the restriction of $\hat x$ to $Q_t$. Then by definition we have that $\hat x^t \in F_t$.  Thus we can set
$\hat \lambda^t_{\hat x^t} = 1$ and $\hat \lambda^t_{v} = 0$ for any other $v \in F_t$, and for any $S \subseteq Q_t$ $\hat Z_S = \zeta^{supp(\hat x^t, Q_t)}_S$. The last equation simply states that $\hat Z_S = 1$ iff $S \subseteq supp(\hat x, V)$, a consistent
definition across $t \in V(T)$.  Hence indeed $(\hat Z, \hat \lambda)$ is
feasible for LPz and  attains $\hat Z_{\{j\}} = \hat x_j$ for each $j \in V$, as
desired.  Note that, effectively, we have argued that the restriction of $\zeta^{supp(\hat x, V)}$ to $\R^{2^{Q_t}}$ for $t \in V(T)$ yields a feasible solution to LPz.

Next we argue that (\ref{consistency0})-(\ref{tconvexity0}) defines an integral
polyhedron.  This is a consequence of the following result, which can
be obtained from Lemma 8.18 of \cite{Laurentima}, although it is not 
stated there in the language of constraints
(\ref{consistency0})-(\ref{tconvexity0}).

\begin{lemma}\label{gifttoMonique}  Suppose that $(Z, \lambda)$ is a feasible solution to
  (\ref{consistency0})-(\ref{tconvexity0}).  Then there exists a vector $W \in \R^{2^{V}}$,
  nonnegative values $\theta_1, \ldots, \theta_k$ and vectors $y^1, \ldots, y^k$ in $\R^{V}_+$ such
  that:
  \begin{itemize}
  \item [(1)] $Z_S = W_S$, for all $S \in \cup_{t \in V(T)}2^{Q_t}$.
  \item[(2)] $\sum_{i = 1}^k \theta_i \ = \ 1$.
  \item[(3)] $y^i$ is feasible for GB, for $1 \le i \le k$.
    \item[(4)] $W \ = \ \sum_{i = 1}^k \theta_i ~\zeta^{supp(y^i, V)}$.
  \end{itemize}
\end{lemma}
As a consequence of (2)-(4), the vector $W$ is a convex combination of the
vectors $\zeta^{supp(y^i, V)}$ which as argued above yield feasible solutions to
LPz, thus yielding the desired result.  We remark that the proof of Lemma 8.18 of \cite{Laurentima} is
related to that of the tree-junction theorem; this technique, evocative of
dynamic programming, was also used in \cite{tw} in a closely related setting.

To complete the proof of Theorem \ref{pbtheorem} we note that the total quantity
of variables and constraints in LPz is $O(\sum_{t \in V(T)} 2^{|Q_t|}) = O(|V(T)| 2^{\omega})$.  This yields part (a) of Theorem \ref{pbtheorem}; part (b) follows by using
standard algorithmic techniques.

\section{Acknowledgments} Many thanks to G\'{a}bor Pataki for 
useful comments.  This research
was partially funded by ONR award N00014-13-1-0042, LANL award ``Grid Science'' and DTRA award  HDTRA1-13-1-0021.
\bibliographystyle{siam}

\bibliography{pwidth}

\newpage

\appendix

\section{Dependence on $\epsilon$ in Theorem \ref{potheorem}} \label{app:epsilon}

In this section we will prove that the two characteristics of Theorem \ref{potheorem} regarding $\epsilon$ (approximation notion and running time) cannot be improved.

First, suppose that there is an algorithm $\A$ such that any PO whose intersection graph has tree-width $\le 2$ can be solved in polynomial time to some given feasibility tolerance $\epsilon < 1$, that is to say the algorithm guarantees $f_i(x) \ge -\epsilon$ for any constraint $f_i(x) \ge 0$. Note that since $\epsilon$ is fixed in this case, the formulation in Theorem \ref{potheorem} yields an algorithm that runs in polynomial time (see the result on Theorem \ref{pbtheorem} for the time it takes to build the LP formulation) but with a weaker
approximation guarantee than the hypothetical algorithm $\A$.

We claim that the existence of algorithm $\A$ implies P = NP.  Consider
the subset-sum problem: given $n\geq 2$ positive integers $a_1, \ldots , a_n$ find
$I \subseteq \{1,\ldots, n\}$ such that $\sum_{j \in I} a_j = \sum_{j \notin I}a_j$. Denoting
$$S \doteq \frac{1}{2} \sum_{j = 1}^n a_j \quad \mbox{and} \quad M \doteq 4n S,$$
the subset-sum problem can be cast as the following (pure feasibility) PO:
\begin{subequations} \label{subsetsum}
  \begin{eqnarray}
    M S y_1 & = & M a_1 x_1, \label{sub1} \\
    M S y_i & = & M a_i x_i \ + \ M S y_{i-1}, \quad 2 \le i \le n,\\
    M y_n & = & M, \label{subn}\\
    M x_j(1 - x_j) & = & 0, \quad 1 \le j \le n, \label{binary}\\
    0 \ \le \ y_j & \le & 1, \quad 1 \le j \le n, \label{ybounded}\\
    0 \ \le \ x_j & \le & 1, \quad 1 \le j \le n. \label{xbounded}
 \end{eqnarray}
\end{subequations}  
\noindent Given a solution $(x, y)$ to \eqref{subsetsum} it is clear that
$x \in \{0,1\}^n$ and that $\sum_{j} a_j x_j = \frac{1}{2} \sum_{j = 1}^n a_j$.
Moreover, the intersection graph of \eqref{subsetsum} has tree-width 2.

By assumption, algorithm $\A$ will produce a solution $(\hat x, \hat y)$ that
violates each of the constraints \eqref{sub1}-\eqref{binary} by at most
$\epsilon$ and that satisfies \eqref{ybounded}-\eqref{xbounded}.   Then adding (\ref{sub1})-(\ref{subn}) yields
\begin{eqnarray}
 \left| \sum_{j = 1}^n a_j \hat x_j  - S \right| &\leq& \frac{\epsilon n}{M}. \label{stupid1}
\end{eqnarray}
Moreover, by \eqref{binary} and \eqref{xbounded} for each $1 \le j \le n$,
$$ \mbox{either $\ 0 \le \hat x_j \le \frac{2 \epsilon}{M} \ $ or $ \ 1 - \frac{2 \epsilon}{M} \le \hat x_j \le 1$.}$$
[This follows from the fact $g(x)= x(1-x)$ is strictly increasing in $[0,1/2)$,
    strictly decreasing in $(1/2,1]$, and $g(2\epsilon/M) > \epsilon/M$.] 
Thus, suppose we round each $\hat x_j$ to the nearest integer, obtaining binary
values $\tilde x_j$ for $1 \le j \le n$.  Using (\ref{stupid1}) we obtain
$$\left| \sum_{j = 1}^n a_j \tilde x_j  -  S \right| \leq \frac{\epsilon n}{M} + \left(\sum_{j = 1}^n a_j \right) \frac{2 \epsilon}{M} \ =  \frac{\epsilon n}{M} + \frac{4S \epsilon}{M}$$
and therefore
$$ \left | \sum_{j = 1}^n a_j \tilde x_j  \, - \,   S \right| \ < \ 1.$$
Since the left hand side of the inequality must be an integer, we conclude that 
$$  \sum_{j = 1}^n a_j \tilde x_j  =   S;$$
which proves that unless $P = NP$ algorithm $\A$ does not exist.

Next, suppose that now that there is an algorithm $\A$ that,
for any $\epsilon < 1$ solves PO problems
to scaled tolerance $\epsilon$ (i.e. the violation of any constraint $f_i(x) \ge 0$ is at most $\epsilon \| f \|_1$) but whose running time
is polynomial, i.e. in particular it depends polynomially
on $\log \epsilon^{-1}$.  This is in contrast with the formulation in Theorem \ref{potheorem} yields an algorithm that runs time polynomial on $n$, $m$ and $\epsilon^{-1}$. Consider an unscaled version of the previous formulation of the subset-sum problem, i.e:
\begin{subequations} \label{subsetsum2}
  \begin{eqnarray}
     S y_1 & = &  a_1 x_1, \label{sub2} \\
     S y_i & = &  a_i x_i \ + \  S y_{i-1}, \quad 2 \le i \le n,\\
     y_n & = & 1, \label{subn2}\\
     x_j(1 - x_j) & = & 0, \quad 1 \le j \le n, \label{binary2}\\
    0 \ \le \ y_j & \le & 1, \quad 1 \le j \le n, \label{ybounded2}\\
    0 \ \le \ x_j & \le & 1, \quad 1 \le j \le n. \label{xbounded2}
 \end{eqnarray}
\end{subequations}  
Define $\epsilon = 1/(3SM)$ and use algorithm $\A$ to find a solution $(\hat x, \hat y)$ that is scaled-$\epsilon$-feasible. Since the
1-norm of any polynomial in constraints \eqref{subsetsum2} is at most $2S+1$, we get that for each constraint $f_i(x,y) \geq 0$
$$ f_i (\hat{x}, \hat{y}) \geq -\epsilon \| f_i \|_1 \geq -\epsilon (2S+1) \geq \frac{1}{M}$$

This way we can reuse the same argument as before to obtain a solution to the subset-sum problem. Since we assume the running time depends on $\log(\epsilon^{-1})$ we get a running time that depends polynomially on $\log(nS)$ yielding the
same contradiction as before.

\section{Proof of part (b) of Theorem \ref{genbtheorem_0}}\label{counting}
Here we describe a procedure that constructs formulation \eqref{LPzeta}
which requires $2^{\omega}m$ oracle queries and with 
additional workload $O^*(2^{\omega}\omega (n + m) + \omega mn)$, as
per Theorem \ref{genbtheorem_0} (b).  Here, as per the formulation,
we have a tree-decomposition $(T,Q)$ of the intersection graph of a problem
$GB$, of width $\omega$.  The critical element in the procedure is the construction of the
sets $F_t$ used in equation \eqref{consistency0}, and we remind the
reader of Definition \ref{bigdef}, and that
for $1 \le i \le m$ constraint $i$ has support $K[i]$ and the set of
feasible solutions for constraint $i$ is indicated by $S^i \subseteq \{0,1\}^{K[i]}$.  Note that $|K[i]| \le \omega$ for all $i$. The procedure operates as follows:

\noindent {\bf 1.} For each constraint $i$, enumerate each partition
of $K[i]$.  Given a partition $(A_1,A_0)$ if the vector $y \in \subseteq \{0,1\}^{K[i]}$ defined by $y_j = k$ if $j \in A_k$ ($k = 0, 1$) is such
that $y \notin S^i$ (i.e., \textit{not} feasible) then we record the
triple $(i, A_1, A_0)$ as a vector of length $|K[i]| + 1$ [with some abuse of notation]. This process
requires $2^{|K[i]|}$ oracle queries.  This sum of all these quantities is
$O(2^{\omega} m)$ but the more precise estimate will be needed.

\noindent {\bf 2.} Let $\cL$ be the list of all vectors recorded in Step 1,
sorted lexicographically; first by the index $i$, then by 
$A_1$ and then by $A_0$.  After postprocessing if necessary, we can
assume that $\cL$ contains no duplicates.  These can be performed in time $O(\omega |\cL| \log |\cL|) = O(2^{\omega} m (\omega + \log m))$.

\noindent {\bf 3.} For each $t \in V(T)$ construct a list of all constraints $i$ such that $K[i] \subseteq Q_t$.  This can be done in time $O(\omega m n)$.

\noindent {\bf 4.} For each $t \in V(T)$ we form the sublist of $\cL$ consisting of all vectors $(i, A_1, A_0)$ (constructed in Step 1) such
that $K[i] \subseteq Q_t$.  Note that for any such $i$ the total number of such vectors is
at most $2^{|K[i]|}$.  Given a vector $(i, A_1, A_0)$ thus enumerated,
we form all vectors of the form $(A'_1, A'_0)$ such that $A'_1 \cup A'_0 = Q_t$ and $A_k \subseteq A'_k$ for $k = 0, 1$.  Let $\cL^t$ be the
list of all vectors obtained this way.  Clearly, $|\cL^t| \le 2^{\omega} |\left\{i \, : \, K[i] \subseteq Q_t \right\}|$.  We lexicographically
sort $\cL^t$.  

\noindent {\bf 5.} For each $t \in V(T)$ we enumerate all vectors $y \in \{0,1\}^{Q_t}$.  For any such  vector $y$, we have that $y \in F_t$ if and only if $y$ is \textit{not} found in the list $\cL^t$; and this test
can be performed in time $O(\omega \log |\cL^t|)$ after lexicographically
sorting the list.\\

The total amount of work entailed in Step 4, using $|Q_t| \le \omega$ for each $t \in V(T)$, is $$O(\sum_{t} \omega |\cL^t| \log |\cL^t|) = O^*( 2^{\omega} \omega \sum_t  |\left\{i \, : \, K[i] \subseteq Q_t \right\}|).$$  Likewise, Step 5 requires $O(\omega 2^{\omega} \sum_t \log |\cL^t|) = O^*(\omega^2 2^{\omega})$.  This completes the proof.

\section{Alternative formulation for Theorem \ref{pbtheorem}} \label{pbproof2}
Here we construct a second formulation that also yields the result of Theorem \ref{pbtheorem}. The number of constraints and
variables in this new construction is upper bounded by those for LPz.

\subsection{Additional Definitions and Second Formulation}
\begin{Def}\label{bigdef2} Let $t \in V(T)$. We let \bmath{\Omega_t} denote the set of pairs $(Y,N)$ with
$Y \cap N = \emptyset, \ Y \cup N \subseteq Q_t,$ and such that 
\begin{enumerate}
\item $|Y| \le 1$  and $|N| = 0$, or
\item $(Y, N)$ partition $Q_t \cap Q_t'$, for some $t' \in V(T)$ with $\{t, t'\} \in E(T)$.
\end{enumerate}
\end{Def}

\noindent \noindent The formulation is as follows. The variables are:
\begin{itemize}
\item A variable $\lambda^t_v$, for each $t \in V(T)$ and each vector $v \in F_t$.
\item A variable $X[Y,N]$, for each pair $(Y, N) \in 2^{n} \times 2^{n}$ with $(Y,N) \in \Omega_t$ for some $t \in V(T)$.
\end{itemize}
\begin{subequations} \label{LPGB}
\begin{eqnarray}
\mbox{(LP-GB):} && \min ~ \sum_{j =1}^n c_j X[\{j\}, \emptyset] \\
\mbox{s.t.} \quad \forall t \in V(T):&& \nonumber \\
&& X[Y, N] \ = \ \sum_{v\in F_t} \lambda^t_v \prod_{j\in Y}  v_j \prod_{j\in N} (1-v_j)\quad  \forall \ (Y,N) \in \Omega_t  \label{consistency}\\
&& \nonumber \\
&& \sum_{v \in F_t} \lambda^t_v \ = \ 1, \quad \lambda^t \ge 0. \label{tconvexity}
\end{eqnarray}
\end{subequations}
We will show below that (a) LP-GB is a relaxation of GB and (b) the relaxation
is exact and that the polyhedron defined by (\ref{consistency})-(\ref{tconvexity}) is integral.
\begin{remark}\label{remmie}
\begin{itemize}
\item[(f.1)] When $(Y,N)$ partition $Q_t \cap Q_{t'}$ for some edge $\{t, t'\}$ then
variable $X[Y,N]$ will appear in the constraint (\ref{consistency}) arising from 
$t$ and also that corresponding to $t'$.  This implies an equation involving
$\lambda^t$ and $\lambda^{t'}$.
\item[(f.2)] The sum on the right-hand side of constraint (\ref{consistency}) could be empty.  This will be the case if
for any $v \in \{0,1\}^{Q_t}$ with $v_j = 1$ for all $j \in Y$ and $v_j = 0$ for all $j \in N$ there exists $1 \le i \le m$ with $K[i] \subseteq Q_t$ and yet
$v_{K[i]} \notin S^i$.  Then (\ref{consistency}) states $X[Y,N] = 0$.
\item[(f.3)] For $v \in F_t$ define
$Y = \{ j \in Q_t \, : \, v_j = 1 \}$ and $N = Q_t - Y$.  Then (\ref{consistency}) states $X[Y,N] = \lambda^t_{v}$.
\item[(f.4)] When $Y = N = \emptyset$ the right-hand side of (\ref{consistency}) 
is $\sum_{v \in F_t} \lambda^t_v$.  Hence we will have $X[\emptyset, \emptyset] = 1$.
\item[(f.5)] The $\lambda$ variables are the same as those in formulation (\ref{LPzeta}). For any edge $\{t,t'\} \in E(T)$ and $Y \subseteq Q_t \cap Q_{t'}$ the terms on the right-hand side of the row (\ref{consistency}) are a subset of the terms on the right-hand side of the row (\ref{consistency0})  corresponding to $Y$ 
(additional statements are possible).
\end{itemize}
\end{remark}
First we show that LP-GB is a relaxation for GB, in a strong sense.
\begin{lemma}\label{basic1}  Let $\tilde x$ be a feasible solution to an instance for GB.
\begin{itemize}
\item [(i)] There is
a feasible, $0/1$-valued solution $(\tilde X, \tilde \lambda)$ to (\ref{LPGB}) such that for each variable $X[Y,N]$ in (\ref{LPGB}) we have $\tilde X[Y, N] = \prod_{j \in Y} \tilde x_j \, \prod_{j \in N} (1 - \tilde x_j)$.
\item [(ii)] As a corollary $ \sum_{j =1}^n c_j \tilde X[\{j\}, \emptyset] \ = \ c^T \tilde x$.
\end{itemize}
\end{lemma}
\noindent {\em Proof.} (i) For each variable $X[Y,N]$ in problem (\ref{LPGB}) we set $\tilde X[Y, N] = \prod_{j \in Y} \tilde x_j ~ \prod_{j \in N} (1 - \tilde x_j)$. Further, for each $t \in V(T)$ let 
$\tilde v(t) \in \{0,1\}^{Q_t}$ be the restriction of $\tilde x$ to $Q_t$, i.e. $\tilde v(t)_j = \tilde x_j$ for each $j \in Q_t$.   
Since $\tilde x$ is feasible, $\tilde v(t) \in F_t$. 
Then we set $\tilde \lambda^t_{\tilde v(t)} = 1$ and $\tilde \lambda^t_{v} = 0$ for every 
vector $v \in F_t$ with $v \neq \tilde v(t)$.
By construction
for every $t \in V(T)$ and $(Y,N) \in \Omega_t$ we have
$\tilde X[Y, N] = 1$  iff $\tilde v(t)_j = 1$ for all $j \in Y$ and 
$\tilde v(t)_j = 0$ for all $j \in N$; in other words (\ref{consistency}) is satisfied. 

\noindent (ii) This follows from (i).
\QED

As a consequence of Lemma \ref{basic1}, Theorem \ref{pbtheorem} will follow
if we can prove that the constraint matrix in (\ref{LPGB}) defines an integral polyhedron.  This will be done in Lemma \ref{goback} given below.  In what follows, we will view $T$ as rooted, i.e. all edges are directed so that $T$ contains a 
directed path from an arbitrarily chosen {\bf leaf} vertex $r$ (the root of $T$) to every other 
vertex.  If $(v, u)$ is an edge thus directed, then we say that $v$ is the parent
of $u$ and $u$ is a child of $v$.  

\begin{Def} A rooted subtree $\tilde T$ is a subtree of $T$, such that  there exists a vertex
$u$ of $\tilde T$ so that $\tilde T$  contains a directed path from $u$ to
every other vertex of $\tilde T$.  We then say that $\tilde T$ is rooted at
$u$.
\end{Def}

\begin{Def}\label{thedefs} Let $\tilde T$ be a rooted subtree of $T$.
\begin{itemize}
\item[(a)] We denote by \bmath{\Omega(\tilde T)} the set $\bigcup_{t \in \tilde T} \Omega_t$.
\item[(b)]We denote by \bmath{\V(\tilde T)} the set $\{ j \ : \, j \in Q_t \ \mbox{for some $t \in \tilde T$}\}$.
\end{itemize}
\end{Def}

\noindent Below we will prove the following result:
\begin{theorem} \label{inductive} Let $(\hat X, \hat \lambda)$ be a feasible solution to the LP-GB problem (\ref{LPGB}). Then for every rooted subtree $\tilde T$ there is a family of vectors
$$ p^{k,\tilde T} \in \{0, 1\}^{\Omega(\tilde T)},$$
vectors
$$x^{k,\tilde T} \in \{0, 1\}^{\V(\tilde T)}$$
and reals
$$ 0 < \mu^{k,\tilde T} \le 1,$$ 
$(k = 1, 2, \ldots, n(\tilde T))$ 
satisfying the following properties:
\begin{itemize}
\item [(a)] For each $1 \le k \le n(\tilde T)$ and each constraint $1 \le i \le m$ of problem GB, if $K[i] \subseteq Q_t$ for some
$t \in \tilde T$, then $x^{k,\tilde T} \in S^i$.
\item [(b)] For $1 \le k \le n(\tilde T)$ and each pair $(Y, N) \in \Omega(\tilde T)$, 
$$ p^{k,\tilde T}[Y, N] \ = \ \prod_{j \in Y} x^{k,\tilde T}_j \ \prod_{j \in N} \left(1 - x^{k,\tilde T}_j \right). $$
As a result, for each $1 \le k \le n(\tilde T)$ and $j \in \V(\tilde T)$, $x^{k,\tilde T}_j = p^{k,\tilde T}[\{j\}, \emptyset]$.
\item [(c)] $\sum_{k = 1}^{n(\tilde T)} \mu^{k,\tilde T} \ = \ 1$.
\item [(d)] For each $(Y,N) \in \Omega(\tilde T)$,
$$ \hat X[Y, N] \ = \ \sum_{k = 1}^{n(\tilde T)} \mu^{k,\tilde T} p^{k,\tilde T}[Y, N].$$
\end{itemize}
\end{theorem}
The family of vectors $p^{k, \tilde T}$ and reals $\mu^{k,\tilde T}$ will be called a 
{\bf decomposition of \bmath{(\hat X, \hat \lambda)} over \bmath{\tilde T}}. \\

\noindent Pending a proof of Theorem \ref{inductive},  we can show that 
the polyhedron defined by the constraints in LP-GB is integral.
\begin{lemma}\label{goback}  The polyhedron defined by (\ref{consistency})-(\ref{tconvexity})
is integral and problems GB and LP-GB have the same value.  
\end{lemma}
\noindent {\em Proof.} Let $(\hat X, \hat \lambda)$ be a feasible solution to LP-GB.  We apply Theorem \ref{inductive} with $\tilde T = T$ obtaining a family of vectors 
$p^{k} \in \{0, 1\}^{\Omega(r)}$, vectors $x^{k} \in \{0,1\}^n$ and reals $\mu^{k}$, for $1 \le k \le n(r)$, satisfying conditions 
(a)-(d) of the theorem.  By (a) and Remark \ref{cliques}, each vector $x^{k}$ is feasible for GB.
By (d), the vector $\hat X$ is a convex combination of the vectors 
$p^{k}$.  This completes the proof, using Remark \ref{remmie} (f.3) to handle
the $\lambda^t$ variables. \QED

This result completes the proof of Theorem \ref{pbtheorem}, pending Theorem \ref{inductive}.

\subsection{Proof of Theorem \ref{inductive}}

Assume
we have a feasible solution $(\hat X, \hat \lambda)$ to (\ref{LPGB}). 
The proof of Theorem \ref{inductive} will be done by induction on the size of $\tilde T$.
First we handle the base case.  

\begin{lemma} \label{basecase}  If $\tilde T$ consists of a single vertex $u$ there is a decomposition of $(\hat X, \hat \lambda)$ over $\tilde T$.
\end{lemma}

\noindent {\em Proof.}  We have that $\Omega(\tilde T) = \Omega_u$ (see Definition \ref{thedefs}).  
By (\ref{tconvexity}) we have $\sum_{v \in F_u} \hat \lambda^u_v = 1$.
Let $n(\tilde T) > 0$ be the number of elements $v \in F_u$ with $\hat \lambda^u_v > 0$ and denote these vectors by
$\{ w(1), \ldots, w(n(\tilde T)) \}$. Then,  
for $1 \le k \le n(\tilde T)$  let $x^{k,\tilde T} = w(k)$ and $\mu^{k,\tilde T} = \hat \lambda^u_{w(k)}$.  
Finally, for $1 \le k \le n(\tilde T)$ we define the vector $p^{k, \tilde T} \in \Omega_u$ by setting
$$ p^{k, \tilde T}[Y,N] \ = \ \prod_{j \in Y} x^{k,\tilde T}_j \ \prod_{j \in N} \left(1 - x^{k,\tilde T}_j \right) $$
for each pair $(Y,N) \in \Omega_u$.  
Now we will verify that conditions (a)-(d) of Theorem \ref{inductive} hold.  Clearly (a)-(c) hold
by construction.  To see that (d) holds, note that $(\hat X, \hat \lambda)$ satisfies (\ref{consistency}), i.e.,
\begin{eqnarray}
 \hat X[Y, N] & = & \sum_{v\in F_u}  \hat \lambda^u_v \prod_{j \in Y} v_j \ \prod_{j \in N} \left(1 - v_j \right)\nonumber \\
 & = & \sum_{k=1}^{n(\tilde{T})}  \mu^{k,\tilde T}\prod_{j \in Y} x^{k,\tilde T}_j \ \prod_{j \in N} \left(1 - x^{k,\tilde T}_j \right)\nonumber \\
& = & \sum_{k = 1}^{n(\tilde T)} \mu^{k,\tilde T} p^{k,\tilde T}[Y, N] \nonumber
\end{eqnarray}
which is condition (e), as desired. \QED\\

Next we prove the general inductive step needed to establish Theorem \ref{inductive}.  The technique used here is related to the junction tree theorem, is similar to one used in \cite{tw} and is
reminiscent of Lemma L of \cite{Laurentima}.\\

Consider a vertex $u$ of $T$ and a subtree $\tilde T$ rooted at $u$ with more than one vertex.    Let $v$ be a child of $u$.  We will apply induction by partitioning $\tilde T$ into two subtrees: 
the subtree $L$ consisting of $v$ and all its descendants in $\tilde T$, and 
the subtree $H = \tilde T - L$.  
Consider a decomposition of $(\hat X, \hat \lambda)$ over $L$ given by the vectors $p^{k,L} \in \{0, 1\}^{\Omega(L)}$
and the positive reals $\mu^{k,L}$  for $k = 1, 2, \ldots, n(L)$, 
and a decomposition of $(\hat X, \hat \lambda)$ over $H$ given by the vectors $p^{k,H} \in \{0, 1\}^{\Omega(H)}$
and the positive reals $\mu^{k,H}$  for $k = 1, 2, \ldots, n(H)$.

Denote by $\p$ the set of partitions of $Q_u \cap Q_v$ into two sets.  Thus, by 
Definition \ref{bigdef2}, for each  
$(\alpha, \beta) \in \p$ we have a variable $X[\alpha, \beta]$.  
Note that $\Omega(\tilde T) = \Omega(H) \cup \Omega(L)$. We construct a family of vectors and reals satisfying
(a)-(d) Theorem \ref{inductive} for $\tilde T$, as follows.

 For each $(\alpha, \beta) \in \p$ such that $\hat X[\alpha, \beta] > 0$, and each pair $i$, $h$ such that
$1 \le i \le n(L)$,  $1 \le h \le n(H)$, and $p^{h, H}[\alpha, \beta] = p^{i, L}[\alpha, \beta] = 1$ we create
a vector $q^{\alpha,\beta}_{ih}$ and a real $\gamma_{ih}^{\alpha, \beta}$ using the rule:\\

\noindent For any vertex $t$ in $\tilde T$ and $(Y, N) \in \Omega_t$: 
\begin{itemize}
\item[(r.1)] If $t \in V(L)$ we set $q^{\alpha,\beta}_{ih}[Y,N] = p^{i, L}[Y,N]$. 
\item[(r.2)] If $t \in V(H)$ we set $q^{\alpha,\beta}_{ih}[Y,N] = p^{h, H}[Y,N]$. 
\end{itemize}
Further, we set
$$ \gamma^{\alpha,\beta}_{ih} \ = \ \frac{ \mu^{i,L} \, \mu^{h, H}} {\hat X[\alpha,\beta]}.$$

\noindent To argue that this construction is valid we note that
since $\hat X[\alpha, \beta] > 0$, pairs of indices $i$, $h$ as listed above must exist, by (d) of
the inductive assumption applied to $H$ and $L$.  Furthermore, we have $\gamma_{ih}^{\alpha, \beta} > 0$.\\

Now we will prove that the $q_{ih}$ and the $\gamma_{ih}$ provide a decomposition of $(\hat X, \hat \lambda)$ over
$\tilde T$.  Let $i$ and $h$ be given.  Since the restriction of $p^{i,L}$ (and $p^{h,H}$) to $L$ (resp., $H$) 
satisfy (a) and (b) of the inductive assumption, so will $q_{ih}$. Thus, there remains to prove (c) and (d).  

First, consider (d). 
Let $(Y, N) \in \Omega(\tilde T)$, say $(Y, N) \in \Omega(H)$.  We claim that
\begin{subequations}
\begin{eqnarray}
\hspace{-30pt} \sum_{\alpha, \beta, i, h} \gamma^{\alpha,\beta}_{ih} q^{\alpha,\beta}_{ih}[Y,N] & = & 
 \sum_{(\alpha, \beta) \in \p \, : \, \hat X[\alpha, \beta] > 0} \sum_{i = 1}^{n(L)} \sum_{h = 1}^{n(H)} \frac{ \mu^{i,L} \, \mu^{h, H}} {\hat X[\alpha,\beta]} p^{i, L}[\alpha, \beta] p^{h, H}[\alpha, \beta] p^{h, H}[Y, N].
\label{intermediate}
\end{eqnarray}
\end{subequations}
This equation holds because in any nonzero term in either expression we must have 
$p^{i, L}[\alpha, \beta] = p^{h, H}[\alpha, \beta] = 1$ and since $(Y,N) \in \Omega(H)$ we also have that
$q^{\alpha,\beta}_{ih}[Y,N]  = p^{h, H}[Y, N]$.  \\
\noindent Now the right-hand side of (\ref{intermediate}) equals
\begin{subequations}
\begin{eqnarray}
\hspace{-30pt} && \sum_{(\alpha, \beta) \in \p \, : \, \hat X[\alpha, \beta] > 0} \left[ \left( \sum_{i = 1}^{n(L)} \frac{\mu^{i,L} p^{i, L}[\alpha, \beta]}{\hat X[\alpha,\beta]} \right) \left( \sum_{h = 1}^{n(H)} \mu^{h,H} p^{h, H}[\alpha, \beta] p^{h, H}[Y, N]\right) \right] \ = \ \\
&& \sum_{(\alpha, \beta) \in \p \, : \, \hat X[\alpha, \beta] > 0} \left( \sum_{h = 1}^{n(H)} \mu^{h,H} p^{h, H}[\alpha, \beta] p^{h, H}[Y, N]\right), \label{pause1}
\end{eqnarray}
\end{subequations}
by the inductive assumption (d) applied to subtree $L$.  The expression in (\ref{pause1}) equals
\begin{subequations}
\begin{eqnarray}
&& \sum_{h = 1}^{n(H)} \left( \left[ \sum_{(\alpha, \beta) \in \p \, : \, \hat X[\alpha, \beta] > 0} p^{h, H}[\alpha, \beta] \right] \mu^{h,H}p^{h, H}[Y, N] \right). \label{ooh}
\end{eqnarray}
\end{subequations}
By inductive property (b) applied to subtree $H$, given $1 \le h \le n(H)$ we have that 
$p^{h,H}[\alpha, \beta] = 1$ for exactly one partition $(\alpha, \beta) \in \p$, and so expression (\ref{ooh}) equals
\begin{eqnarray}
&& \sum_{h = 1}^{n(H)}  \mu^{h,H}p^{h, H}[Y, N]. \label{aah}
\end{eqnarray}
In summary, 
$$ \sum_{\alpha, \beta, i, h} \gamma^{\alpha,\beta}_{ih} q^{\alpha,\beta}_{ih}[Y,N] \ = \ \sum_{h = 1}^{n(H)}  \mu^{h,H}p^{h, H}[Y, N]. $$
and by induction applied to the subtree $H$ this quantity equals $\hat X[Y, N]$.   Thus property (d) does indeed hold. 

Finally we turn to (c).  Inductively, (c) and (d) hold for trees $L$ and $H$.  
Thus, as noted in Remark \ref{remmie}(f.4), $X[\emptyset, \emptyset] = 1$.  But
we have just shown that (d) holds for $\tilde T$, and in particular that it holds
for $Y = N = \emptyset$.  Using Remark \ref{remmie}(f.4) we obtain that (c) holds
for $\tilde T$, as desired.
\QED

\end{document}